\documentclass[11pt]{amsart}
\usepackage{tocvsec2}
\usepackage{amsmath,amsfonts,latexsym,amssymb,amsthm,mathrsfs,array}

\usepackage[letterpaper]{geometry}
\usepackage[usenames]{color}
\usepackage[colorlinks=true,linkcolor=webgreen,filecolor=webbrown,citecolor=webgreen]{hyperref}

\definecolor{webgreen}{rgb}{0,.5,0}
\definecolor{webbrown}{rgb}{.6,0,0}

\newcommand{\CC}{{\mathbb C}}
\newcommand{\ZZ}{{\mathbb Z}}

\def\CC{{\mathbb C}}
\def\ZZ{{\mathbb Z}}

\newtheorem{dfn}{Definition}[section]
\newcommand{\bdfn}{\begin{dfn}\rm}
\newcommand{\edfn}{\end{dfn}}
\newtheorem{thm}[dfn]{Theorem}
\newcommand{\bthm}{\begin{thm}}
\newcommand{\ethm}{\end{thm}}
\newtheorem{lmma}[dfn]{Lemma}                   
\newcommand{\blmma}{\begin{lmma}}                   
\newcommand{\elmma}{\end{lmma}}                   
\newtheorem{ppsn}[dfn]{Proposition}
\newcommand{\bppsn}{\begin{ppsn}}
\newcommand{\eppsn}{\end{ppsn}}
\newtheorem{crlre}[dfn]{Corollary}
\newcommand{\bcrlre}{\begin{crlre}} 
\newcommand{\ecrlre}{\end{crlre}}
\newtheorem{rmk}[dfn]{Remark}
\newcommand{\brmk}{\begin{rmk}\rm} 
\newcommand{\ermk}{\end{rmk}}

\numberwithin{equation}{section}

\title[Classification of irreducible Harish-Chandra modules]
{Classification of irreducible Harish-Chandra modules
over full toroidal Lie algebras and higher-dimensional Virasoro algebras}
\author{Souvik Pal}         
\address{Souvik Pal, Department of Mathematics, Indian Institute of Science, CV Raman Road, Bengaluru, Karnataka 560012; and
Department of Sciences and Humanities, CHRIST University, Mysore Road, Bangalore 560 074, India.}     
\email{pal.souvik90@gmail.com, souvik.pal@christuniversity.in}  
\date{}

\begin{document}



\maketitle
\begin{abstract}
In this paper, we classify the irreducible Harish-Chandra
modules over the full toroidal Lie algebra, which is a natural
higher-dimensional analogue of the affine-Virasoro algebra. In particular, we complete the classification of irreducible bounded modules, which were studied by Billig for non-zero level modules [\textit{Int.\ Math.\ Res.\ Not.}\ 2006]. As a by-product, we also obtain the classification of irreducible Harish-Chandra modules over the higher-dimensional Virasoro algebra, which was introduced by Rao--Moody [\textit{Comm.\ Math.\ Phys.}\ 1994], thereby generalizing the well-known result of O. Mathieu [\textit{Invent.\ Math.}\ 1992] for the classical Virasoro algebra. More precisely, we show that any irreducible Harish-Chandra module over the higher-dimensional Virasoro algebra turns out to be either a quotient of a module of tensor fields on a torus or a highest weight type module up to a twist of an automorphism, as conjectured by Eswara Rao in 2004.
\end{abstract} 

\settocdepth{section}
\tableofcontents
\section{Introduction}

The affine Lie algebras have been a great success story. Not only did they effectively generalize the well-developed theory of
finite-dimensional simple Lie algebras to the infinite-dimensional
setting, but these Lie algebras have surprisingly found several
remarkable connections in many diverse branches of mathematics and
physics, like modular forms, vertex algebras, to name a few \cite{FK,K,KRR}.
An affine Lie algebra can be realized as the (1-dimensional) universal
central extension of algebra of functions on the unit circle $\mathbb{S}^1$ with values
in a finite-dimensional simple Lie algebra. There is another Lie algebra
which is playing an increasingly important role in mathematics and theoretical physics
called the Virasoro algebra. This Lie algebra can be viewed as a
(1-dimensional) central extension of the polynomial vector fields on
$\mathbb{S}^1$ or alternatively (as a 1-dimensional central extension) of
the Lie algebra of derivations of the Laurent polynomial algebra in a
single variable. 

Most applications of affine and Virasoro algebras arise from their
representation theory. The Virasoro algebra also plays a key role in the
representation theory of affine Lie algebras, as it acts on every (except
when the level is negative of the dual Coxeter number) highest weight
module for the affine Lie algebra. This remarkable connection eventually led to constructing the
affine-Virasoro algebra \cite{V,G}, which is the semi-direct product of
the affine Lie algebra and the Virasoro algebra with a common extension.
This Lie algebra has since emerged to be an extremely important
object of study and its connection to conformal field theory has been
explained in great detail in \cite{FMS}. In particular, the even part of
the $N=3$ superconformal algebra is simply the affine-Virasoro
algebra whose underlying finite-dimensional simple Lie
algebra is $\mathfrak{sl}_2$. 

The most important and widely studied class of representations are probably the ones admitting weight space decompositions which separate the modules into finite-dimensional weight spaces, thereby enabling us to compute characters. These representations are usually referred to as Harish-Chandra modules in the literature. The problem of classifying irreducible Harish-Chandra modules for the affine Kac-Moody algebra and the Virasoro algebra was initiated in \cite{C} and \cite{CP} respectively. For recent progress related to the classification of irreducible Harish-Chandra modules over the affine Kac-Moody algebra, see \cite{CP1, CP2, DG, FT, P, E3} and the references therein. The classification of irreducible Harish-Chandra modules over the affine-Virasoro algebra has been recently settled in \cite{CLW,LPX} (also see \cite{GHL}). It is thus natural to generalize the notion of affine-Virasoro algebra in the higher-dimensional set-up and then develop the corresponding theory of Harish-Chandra modules, which is the underlying theme of this paper. 

To generalize the construction of the affine-Virasoro algebra, we first need to find a suitable higher-dimensional analogue of the affine Lie algebra. This is taken care of by the notion of toroidal Lie algebra \cite{EMY}, which is a natural higher-dimensional generalization of the affine Lie algebra. We begin with a finite-dimensional simple Lie algebra $\mathfrak{g}$ and then consider $\mathfrak{g}$-valued polynomial maps from an $(n+1)$-dimensional torus, which are identified with the multiloop algebra $L(\mathfrak{g}) = \mathfrak{g} \otimes A$, where $A = \CC[t_0^{\pm 1}, \ldots, t_n^{\pm 1}]$ ($n \geqslant 1)$. Parallel to the construction of affine Lie algebras, we thereby consider the universal central extension of $L(\mathfrak{g})$ given by $\overline{L}(\mathfrak{g})=L(\mathfrak{g}) \oplus \mathcal{Z}$ (see (\ref{Universal})) to obtain the toroidal Lie algebra. But unlike in the affine case (i.e. when $n =0$), $\mathcal{Z}$ is infinite-dimensional for $n \geqslant 1$. The classification of irreducible Harish-Chandra modules for $\overline{L}(\mathfrak{g})$ can be found in \cite{ML}.

We next provide a suitable higher-dimensional analogue of the
Virasoro algebra. Consider the derivation algebra of $A$,
which we denote by $\mathcal{W}_{n+1}$ (see (\ref{Witt})). In 
contrast to the one variable case, $\mathcal{W}_{n+1}$ is always
centrally closed for $n \geqslant 1$ \cite{RSS}. This unexpected phenomenon
has led to numerous generalizations of the Virasoro algebra, whose
representations have been extensively studied. These include, but are not
limited to, the higher rank Virasoro algebras \cite{LZ1,S}, the
generalized Virasoro algebras \cite{BZ,GLZ1}, the Virasoro-like algebra
\cite{GL,LT} and the $\mathbb{Q}$-Virasoro algebra \cite{M}. Meanwhile
Eswara Rao and Moody constructed the first interesting lowest energy
representation of a non-central extension of $\mathcal{W}_{n+1}$
\cite{EM}. A major development in \cite{EM} is the appearance of a
$\mathcal{Z}$-valued $2$-cocyle $\phi_1$ on $\mathcal{W}_{n+1}$, which
generalizes the Virasoro cocycle. This idea was taken one step
further by Larrson, who constructed a
bigger class of representations for the toroidal Lie algebras in \cite{TL}, by means of a linear combination of $2$-cocycles $\phi_1$ and $\phi_2$, as opposed
to a single $2$-cocycle used in \cite{EM}. We shall use this
linear combination of $\phi_1$ and $\phi_2$ (see (\ref{2-Cocycle}))
to obtain our higher-dimensional analogue of the Virasoro
algebra. Let us denote this particular generalization of the 
Virasoro algebra, via an (infinite-dimensional) abelian extension, by
$\mathcal{V}ir = \mathcal{W}_{n+1} \oplus \mathcal{Z}$ (see Subsection
\ref{Virasoro} and Remark \ref{Generalization}). Note that if
$n=0$, this abelian extension turns out to be a central extension. 

Having constructed $\overline{L}(\mathfrak{g})$ and
$\mathcal{V}ir$, we take their semi-direct product to obtain the full
toroidal Lie algebra $\tau$ (see (\ref{Full Toroidal})), which is clearly
a higher-dimensional generalization of the affine-Virasoro algebra.
Note that in general, $\mathcal{W}_{n+1}$ is not a subalgebra of
$\mathcal{V}ir$. In fact, $\mathcal{W}_{n+1}$ is a ($centerless$)
subalgebra of $\mathcal{V}ir$ precisely when the linear combination
$\phi$ of the $2$-cocycles $\phi_1$ and $\phi_2$ used in our construction
of $\mathcal{V}ir$ is taken to be $0$. On the other hand, it is easy to
see that the centers of both $\mathcal{V}ir$ and $\tau$ are spanned by
exactly the same $(n+1)$ elements $K_0, \ldots, K_n$ (see (\ref{Full
Toroidal}) and (\ref{Higher-dimensional})). If all these $(n+1)$ elements
act trivially on a $\tau$-module (or a $\mathcal{V}ir$-module), we say
that the representation has \textit{level zero}, otherwise we say that it has
\textit{non-zero level}.

The three distinct Lie algebras discussed above, namely $\mathcal{V}ir, \
\tau$ and $\mathcal{W}_{n+1}$ (corresponding to $\phi= 0$) are all
 $\mathbb{Z}^{n+1}$-graded and so a natural class of
Harish-Chandra modules over these Lie algebras are the weight
modules having finite-dimensional weight spaces with respect to this
grading (see Subsection \ref{Representations}). Substantial efforts were
required to classify the irreducible Harish-Chandra modules over
$\mathcal{W}_{n+1}$ (see \cite{Y,BF,LZ1,MZ,E4,E2,S} and the references therein), which eventually resulted in
the final classification, due to a fairly recent paper by
Billig--Futorny \cite{BF}. However, the classification problem for
irreducible Harish-Chandra modules over $\mathcal{V}ir$ is still
unresolved, which has also left the representation theory of full toroidal Lie
algebras somewhat incomplete. The irreducible Harish-Chandra modules for
the classical Virasoro algebra (obtained by putting $n=0$) were
classified by O.\ Mathieu in \cite{OM}, thereby proving a conjecture of
Kac. More precisely, it was shown in \cite{OM} (also see \cite{MP}) that
every non-trivial irreducible Harish-Chandra module over the classical Virasoro algebra
falls into two classes:
(i)~quotients of modules of tensor fields on $\mathbb{S}^1$ or
(ii)~highest/lowest weight modules.
Unlike the classical Virasoro case, there is no standard
triangular decomposition of $\mathcal{V}ir$, as the lattice
$\mathbb{Z}^{n+1}$ ($n \geqslant 1$) does not naturally decompose into
positive and negative parts. Therefore there are no obvious analogues of
highest weight modules with finite-dimensional weight spaces and
so the correct formulation of Kac's conjecture was initially not clear for $\mathcal{V}ir$-modules.

For the first time in \cite{BB} (also see \cite{BZ}), a large family of
highest weight type modules with finite-dimensional weight spaces
was constructed for some classes of $\mathbb{Z}^{n+1}$-graded Lie
algebras, thereby generalizing the notion of highest weight modules over
$\mathcal{V}ir$ (and over $\tau$). 

In \cite{E}, a couple of decades ago, Eswara Rao conjectured that every non-trivial irreducible Harish-Chandra module over $\mathcal{V}ir$ has to come from  either (i) a quotient of a module of tensor fields on an $(n+1)$-dimensional torus or (ii) a highest weight type module twisted by an automorphism from $GL(n+1, \mathbb{Z})$.

The above stated conjecture of Eswara Rao is also intimately connected to the classification of irreducible Harish-Chandra modules for the full toroidal Lie algebra $\tau$. The classification of irreducible \textit{integrable} Harish-Chandra modules over $\tau$ was carried out in \cite{EJ} (also see \cite{CD}). In \cite{B}, several irreducible modules over $\tau$ of highest weight type (referred to as \textit{bounded} modules in \cite{B}) of \textit{non-zero level}  were explicitly constructed through vertex algebras. 
These bounded modules were also used to obtain realizations of highest weight type modules for $\mathcal{W}_{n+1}$ \cite{YF}. 

The main purpose of our paper is to classify the irreducible Harish-Chandra modules over $\tau$. In this process, we also obtain the classification of all possible irreducible Harish-Chandra modules over $\mathcal{V}ir$, thereby proving Eswara Rao's conjecture.        

In addition, we expect our methods to be widely applicable,
especially in the context of extended affine Lie algebras (EALAs). The EALAs initially appeared in the set-up of elliptic singularities \cite{SY} and then in the framework of Lie algebras related to quantum field gauge theory \cite{KT}. 
Owing to remarkable breakthroughs in
\cite{ABFP} and \cite{Ne}, it is now evident that \textit{almost
every} EALA can be in fact realized as an extension (coming from
$\mathcal{W}_{n+1}$) of some central extension of a multiloop algebra
(see toroidal EALA \cite{YB}, Hamiltonian EALA \cite{ER}, contact EALA
\cite[Section 8]{ER} and minimal EALA \cite[Remark 8.3]{ER} for examples
of such EALAs). Moreover the subalgebra of $\mathcal{V}ir$ given by the
infinite-dimensional abelian extension of divergence zero vector fields,
namely $\mathcal{S}_{n+1} \oplus \mathcal{Z}$, where $\mathcal{S}_{n+1} =
\text{span} \{\sum_{i=0}^{n}u_it^{\underline{r}}d_i \ | \
\sum_{i=0}^{n}u_ir_i = 0, \ \underline{u} \in \mathbb{C}^{n+1}, \
\underline{r} \in \mathbb{Z}^{n+1} \}$, is prominently associated with the magnetic hydrodynamic equations \cite{B2}. The
representation theory of EALAs is still in progress and most of the
work done concerns only \textit{integrable} modules (see
\cite{ER} and the references therein) and \textit{non-zero level bounded}
modules \cite{YB}. We hope that this paper will also contribute in the
further development of the representation theory of EALAs, which will give rise to new applications of this interesting class of Lie algebras. In a recent paper, irreducible highest weight modules over the minimal EALA were studied by Shrawan Kumar \cite{Ku} in connection to modular representation theory, while addressing some questions raised by Lusztig.            

\smallskip

\noindent \textbf{Organization of the paper.} After constructing the full toroidal Lie algebra $\tau$ and the higher-dimensional Virasoro algebra $\mathcal{V}ir$ in Section
\ref{Notation}, we classify the irreducible cuspidal modules (weight modules with uniformly bounded weight spaces) in Section \ref{Central}. Unlike the affine-Virasoro (or classical Virasoro) case, the operators in $\mathcal{Z}$ do not act by scalars, but rather as linear operators on these modules. In Theorem \ref{T2.13}, we show that $\mathcal{Z}$ acts trivially on any irreducible cuspidal module over $\tau$ and over $\mathcal{V}ir$. We stress that in the set-up of $\mathcal{V}ir$-modules, the absence of $L(\mathfrak{g})$ poses a much stiffer challenge and so we had to devise new techniques to prove this assertion (see Lemma \ref{L2.12}). We then utilize the theory of $A$-cover to classify the irreducible cuspidal modules over $\tau$ (see Theorem \ref{T2.14}). To deal with the presence of $L(\mathfrak{g})$ inside $\tau$, we construct a new family of differentiators to show that the $A$-cover is also cuspidal in this framework (see Lemma \ref{Cuspidal Cover}).

In Section \ref{GHW modules}, we prove a key lemma (see Lemma \ref{L3.2})
which helps us to show that every irreducible Harish-Chandra module over
$\tau$ is either a cuspidal module or a generalized highest weight (GHW)
module (see Theorem \ref{T2.9}). In Section \ref{Highest Weight}, we define
highest weight type modules and explicitly describe their highest weight
spaces (see Theorem \ref{L(X)}). Contrary to the one variable case, the
(non-trivial) highest weight space in our setting is always
infinite-dimensional.

In Section \ref{S6}, we prove that every irreducible GHW module has to
come from a highest weight type module (see Theorem \ref{T4.2}).
This passage from GHW modules to highest weight type modules is quite
challenging, unlike in the one variable set-up, where these two classes
of modules coincide much more easily. Furthermore one of the main obstacles to earlier approaches with regard to this transition from GHW modules to highest
weight type modules is that the  Witt algebra
$\mathcal{W}_{n+1}$ of rank $(n+1)$ and the solenoidal Lie algebra
$W(\underline{\gamma})$ (see Subsection \ref{Solenoidal}) are not
subalgebras of $\mathcal{V}ir$ and both of them in fact generate the
\textit{infinite-dimensional} Lie algebra $\mathcal{Z}$. Therefore we had to
resort to different methods to prove this result (for instance, see Lemma
\ref{L4.8} and Theorem \ref{T4.10}).

In Section \ref{Theorems}, we provide the classification of
irreducible Harish-Chandra modules over both $\tau$ (see Theorem
\ref{C1}) and $\mathcal{V}ir$ (see Theorem \ref{C2}). As a corollary, we obtain all possible irreducible modules in the category of \textit{bounded} modules, which was introduced in \cite{B} (see Remark \ref{Bounded}). Finally we also classify the irreducible Harish-Chandra modules over $\mathcal{H}\mathcal{V}ir$ (see Remark \ref{Analogue}), which is a natural higher-dimensional analogue of the twisted Heisenberg--Virasoro algebra admitting a common extension. The irreducible Harish-Chandra modules over HVir were classified in \cite{LZ2}. 

\section{Notations and Preliminaries}\label{Notation}

Throughout this paper, all vector spaces, algebras and tensor products are over the field of complex numbers $\mathbb{C}$. We shall denote the set of integers, natural numbers and non-negative integers by $\mathbb{Z}$, $\mathbb{N}$ and $\mathbb{Z_{+}}$ respectively. For any $n \in \mathbb{Z_{+}}$, $\{e_0, \ldots, e_n \}$ is the canonical $\mathbb{Z}$-basis of $\mathbb{Z}^{n+1}$ and for any Lie algebra $L$, its universal enveloping algebra will be denoted by $U(L)$.  
\subsection{Full Toroidal Lie Algebra} \label{Toroidal}
Consider a finite-dimensional simple Lie algebra $\mathfrak{g}$ equipped with a Cartan subalgebra $\mathfrak{h}$. Then $\mathfrak{g}$ is endowed with a symmetric, non-degenerate and associative bilinear form, which is unique up to scalars. We shall denote this bilinear form by $(\cdot|\cdot)$. 

Let $A = \CC[t_0^{\pm 1}, \ldots, t_n^{\pm 1}]$ ($n \geqslant 1)$ be the algebra of Laurent polynomials in $(n+1)$ variables.
Consider the multiloop algebra given by
\[ {
	L(\mathfrak{g})=\mathfrak{g} \otimes A , \qquad
	[x \otimes f, y \otimes g] = [x,y] \otimes fg \ \forall \ x,y \in \mathfrak{g} \ \text{and} \ f, g \in A.
}
\]
For any $x \in \mathfrak{g}$ and $\underline{k} \in \mathbb{Z}^{n+1}$,
write $t^{\underline{k}}=t_0^{k_0} \ldots t_n^{k_n}$ and let $x \otimes
t^{\underline{k}}$ denote a typical element of $L(\mathfrak{g})$. Now
consider the module of differentials ($\Omega_A,d$) of $A$, which is the
free left $A$-module with basis $\{K_0, \ldots, K_n \}$ along with the
differential map $d : A \longrightarrow \Omega_A$. The image of this
map is spanned by
$d(t^{\underline{k}}) = \sum_{i =0}^{n}k_i t^{\underline{k}}K_i$ for
$\underline{k} \in \ZZ^{n+1}$, where we have $K_i = t_i^{-1}dt_i \
\forall \ 0 \leqslant i \leqslant n$. More precisely,
\begin{align*}
	\Omega_A = \text{span} \{t^{\underline{k}}K_i \ | \ 0 \leqslant i \leqslant n, \ \underline{k} \in \ZZ^{n+1} \},\ 
	dA =  \text{span} \big \{\sum_{i=0}^{n}{k_i t^{\underline{k}}K_i} \ | \ \underline{k} \in \ZZ^{n+1} \big \}.
\end{align*}
If we now consider the quotient space $\mathcal{Z} = \Omega_A / dA$, then we know that
\begin{align}\label{Universal}
	\overline{L}(\mathfrak{g})=L(\mathfrak{g}) \oplus \mathcal{Z}
\end{align} 
is the universal central extension of $L(\mathfrak{g})$ \cite{Ka,EMY}. By
abuse of notation, we shall denote the image of $t^{\underline{k}}K_i$ in
$\mathcal{Z}$ again by itself and define the bracket operation
on $\overline{L}(\mathfrak{g})$ as follows:
\begin{enumerate}
	\item $[x \otimes t^{\underline{k}}, y \otimes t^{\underline{l}}] = [x,y] \otimes t^{\underline{k} + \underline{l}} + (x|y) \sum_{i=0}^{n} k_it^{\underline{k} + \underline{l}}K_i$,
	\item $\mathcal{Z}$ is central in $\overline{L}(\mathfrak{g})$.  
\end{enumerate}
The final ingredient needed in defining the full toroidal Lie algebra is
the derivation algebra of $A$ which we denote by
$\mathcal{W}_{n+1}$. Setting $d_i=
t_i\frac{\partial}{\partial{t_i}}$ which acts on $A$ as derivations,  we can also define
\begin{align} \label{Witt}
	\mathcal{W}_{n+1}= \text{span} \{t^{\underline{r}}d_i : \underline{r} \in \mathbb{Z}^{n+1}, \ 0 \leqslant i \leqslant n \}.    
\end{align}
This infinite-dimensional derivation algebra is itself an extremely classical object in its own right and is popularly known as the \textit{Witt algebra of rank} $(n+1)$. It is now easy to verify that
\begin{align*}
	[t^{\underline{r}}d_i, t^{\underline{s}}d_j]=s_it^{\underline{r} + \underline{s}}d_j - r_jt^{\underline{r} + \underline{s}}d_i.
\end{align*}
\brmk \label{R1.1}
If we consider
\[
	F_{i,j} = t_it_j^{-1}d_j, \qquad
	F_{i,n+1} = -t_i \sum_{k=0}^{n} d_k, \qquad
	F_{n+1,i} = t_i^{-1}d_i, \qquad
	F_{n+1,n+1} = - \sum_{k=0}^{n}d_k
\]
for $0 \leqslant i,j \leqslant n$,
then $\mathcal{L} = \{F_{i,j} \ | \ 0 \leqslant i \neq j \leqslant n+1 \} \bigcup \{F_{i,i} - F_{i+1,i+1} \ | \ 0 \leqslant i \leqslant n \}$ forms the standard basis of $\mathfrak{sl}_{n+2}(\mathbb{C})$ inside $\mathcal{W}_{n+1}$. The Cartan subalgebra of this Lie algebra  spanned by $\mathcal{L}$ is given by $D = \text{span} \{d_0, \ldots, d_n\}$, whence its root lattice can be readily identified with $\mathbb{Z}^{n+1}$. 
\ermk
Again any $d \in \mathcal{W}_{n+1}$ can be extended to a derivation on $L(\mathfrak{g})$ by setting
\begin{align*}
	d(x \otimes f)=x \otimes df \ \forall \ x \in \mathfrak{g}, \ f \in A,
\end{align*}
which subsequently has a unique extension to $\overline{L}(\mathfrak{g})$ via the action  
\begin{align*}
	t^{\underline{r}}d_i(t^{\underline{s}}K_j)=s_it^{\underline{r} + \underline{s}}K_j + \delta_{ij} \sum_{p=0}^{n}r_p t^{\underline{r} + \underline{s}}K_p.
\end{align*}
Moreover it is well-known that $\mathcal{W}_{n+1}$ admits two non-trivial
2-cocycles $\phi_1$ and $\phi_2$ with values in $\mathcal{Z}$:
\begin{equation} \label{2-Cocycle}
	\phi_1(t^{\underline{r}}d_i, t^{\underline{s}}d_j) = -s_ir_j
	\sum_{p=0}^{n}{r_p t^{\underline{r} + \underline{s}}K_p}, \qquad
	\phi_2(t^{\underline{r}}d_i, t^{\underline{s}}d_j) = r_is_j \sum_{p=0}^{n}{r_p t^{\underline{r} + \underline{s}}K_p}          
\end{equation}
(see \cite{B} for more details).
Let $\phi$ be any arbitrary linear combination of $\phi_1$ and $\phi_2$. Then we can define the $full$ $toroidal$ $Lie$ $algebra$ $of$ $rank$ $(n+1)$ (relative to $\mathfrak{g}$ and $\phi$) by setting
\begin{align}\label{Full Toroidal}
	\tau := \tau(\phi) = L(\mathfrak{g}) \oplus \mathcal{Z} \oplus \mathcal{W}_{n+1}
\end{align}
with the following bracket operations besides the relations (1) and
(2) found in Subsection \ref{Toroidal}:
\begin{enumerate}
	\item $[t^{\underline{r}}d_i,t^{\underline{s}}K_j]=s_it^{\underline{r} + \underline{s}}K_j + \delta_{ij} \sum_{p=0}^{n}r_p t^{\underline{r} + \underline{s}}K_p$,
	\item $[t^{\underline{r}}d_i,t^{\underline{s}}d_j]=s_it^{\underline{r} + \underline{s}}d_j-r_j t^{\underline{r} + \underline{s}}d_i + \phi(t^{\underline{r}}d_i,t^{\underline{s}}d_j)$,
	\item $[t^{\underline{r}}d_i,x \otimes t^{\underline{s}}]=s_ix \otimes t^{\underline{r} + \underline{s}} \ \forall \ x \in \mathfrak{g},\ \underline{r},\ \underline{s} \in \mathbb{Z}^{n+1}, \ 0 \leqslant i,j \leqslant n$.   
\end{enumerate}
Observe that  $\tau$ is naturally $\mathbb{Z}^{n+1}$-graded, with dim $\tau_{\underline{m}} = \mathrm{dim} \mathfrak{g} + (2n+1)$ for all non-zero $\underline{m} \in \mathbb{Z}^{n+1}$ and dim $\tau_{\underline{0}} = \mathrm{dim} \mathfrak{g} + 2(n+1)$ . Also note that the $center$ of $\tau$ is spanned by $K_0, \ldots, K_n$. 

\subsection{Higher-dimensional Virasoro Algebra} \label{Virasoro}
We have already seen that the construction of the full toroidal Lie algebra $\tau$ in $(n+1)$ variables involves the abelian extension of the Witt algebra $\mathcal{W}_{n+1}$ of rank $(n+1)$. We shall denote this abelian extension by 
\begin{align}\label{Higher-dimensional}
	\mathcal{V}ir := \mathcal{V}ir(\phi) = \mathcal{Z} \oplus \mathcal{W}_{n+1},
\end{align} 
which is clearly a subalgebra of $\tau$ and so the bracket operations on
$\mathcal{V}ir$ are simply induced from that of $\tau$. It is clear that
the $center$ of $\mathcal{V}ir$ is spanned by the $(n+1)$ elements $K_0,
\ldots, K_n$.

\brmk \label{Generalization}
When $n=0$, $\mathcal{V}ir$ is in fact a \textit{central
extension} of the Witt algebra of rank $1$, which can be also viewed as
the Lie algebra of polynomial vector fields on the unit circle.
Consequently $\mathcal{V}ir$ becomes isomorphic to the classical
$Virasoro$ $algebra$. This reveals that $\mathcal{V}ir$ can be regarded
as a natural higher-dimensional analogue of the classical Virasoro
algebra, in the sense that it is an abelian extension of the Lie algebra
of polynomial vector fields on the $(n+1)$-dimensional torus.
So it makes sense to refer to the Lie algebra $\mathcal{V}ir$ as
the $higher$-$dimensional$ $Virasoro$ $algebra$.
\ermk

\subsection{Solenoidal Lie Algebra}\label{Solenoidal}

A vector $\underline{\gamma} = (\gamma_0, \ldots, \gamma_n) \in
\mathbb{C}^{n+1}$ is called \textit{generic} if
$\sum_{i=0}^{n}\gamma_i r_i \neq 0$ for each $\underline{0} \neq
\underline{r} \in \mathbb{Z}^{n+1}$. Consider any generic vector
$\underline{\gamma} \in \mathbb{C}^{n+1}$ and set $D(\underline{\gamma})
= \sum_{i=0}^{n} \gamma_id_i$. Then a ($centerless$) $solenoidal$ $Lie$
$algebra$ $W(\underline{\gamma})$ is the subalgebra of
$\mathcal{W}_{n+1}$, which is defined as
\begin{align*}
	W(\underline{\gamma}) = AD(\underline{\gamma}).	
\end{align*} 
The Lie bracket on $W(\underline{\gamma})$ is induced from that of $\mathcal{W}_{n+1}$ and is given by 
\begin{align*}
	[t^{\underline{r}}D(\underline{\gamma}), t^{\underline{s}}D(\underline{\gamma})] = \sum_{i=0}^{n}(s_i - r_i) \gamma_i t^{\underline{r} + \underline{s}}D(\underline{\gamma}) \ \forall \ \underline{r}, \underline{s} \in \mathbb{Z}^{n+1}. 
\end{align*}

\brmk 
In general  $W(\underline{\gamma})$ need not
be a subalgebra of $\tau$, as the commutator relations involving
$\mathcal{W}_{n+1}$ in $\tau$ also contain elements of the abelian
extension $\mathcal{Z}$. In fact, $W(\underline{\gamma})$ is a subalgebra
of $\tau$ precisely when the linear combination $\phi$ of the $2$-cocycles
$\phi_1$ and $\phi_2$ used in the construction of $\tau$ is $zero$.
Similarly the Witt algebra $\mathcal{W}_{n+1}$ of rank $(n+1)$ is not
necessarily a subalgebra of $\tau$. 
\ermk

\subsection{Twisted Heisenberg--Virasoro Algebra}\label{Twisted}
This is the Lie algebra HVir with basis
\begin{align*}
	\{x_i, I(j), C_D, C_{DI}, C_I \ | \ i, j \in \mathbb{Z} \}
\end{align*}
and the Lie bracket given by
\begin{align*}
	[x_i,x_j] = (j-i)x_{i+j} + \delta_{i+j,0} \dfrac{i^3 - i}{12}C_D, \\
	[x_i, I(j)] = jI(i+j) + \delta_{i+j,0}(i^2 + i)C_{DI}, \\
	[I(i), I(j)] = i \delta_{i+j,0}C_I, \\
	[\text{HVir}, C_D] = [\text{HVir}, C_{DI}] = [\text{HVir}, C_I] = 0. 
\end{align*}

\subsection{Change of Co-ordinates}\label{Co-ordinates}

Let $G=GL(n+1,\mathbb{Z})$ stand for the general linear group of
invertible $(n+1) \times (n+1)$ matrices with integer entries. Then $G$
acts on $\mathbb{Z}^{n+1}$ by means of matrix multiplication. Let us now
fix any $A = (a_{ij})_{0 \leqslant i,j \leqslant n}$ in $G$ and define 
\[
T_{A}(x \otimes t^{\underline{m}}) = x \otimes t^{\underline{m}A^t},
\qquad
T_{A}(t^{\underline{m}}K_j) = \sum_{p=0}^{n}a_{pj}
t^{\underline{m}A^t}K_p, \qquad
T_{A}(t^{\underline{m}}d_j) = \sum_{p=0}^{n}b_{jp}
t^{\underline{m}A^t}d_p,
\]
where $0 \leqslant j \leqslant n$,
$B = (b_{ij}) = A^{-1}, \ A^{t}$ is
the transpose of $A$ and $\underline{m}$ is a \textit{row
vector} in $\mathbb{Z}^{n+1}$. One checks that $T_{A}$ is
an automorphism of $\tau$ that takes our original full
toroidal Lie algebra $\tau$ to a newly formed full toroidal Lie
algebra isomorphic to $\tau$.   This phenomenon is
termed as change of co-ordinates, which enables us to consider the action
of $G$ on $\tau$ via automorphisms. In this paper, we shall use this
notion without further comments and simply refer to it as
\textit{a change of co-ordinates}.

\subsection{Harish-Chandra Modules}\label{Representations}
$V$ is said to be a Harish-Chandra module over $\tau$ if it
satisfies:
\begin{enumerate}
	\item $V = \bigoplus_{\underline{m} \in \mathbb{Z}^{n+1}} V_{\underline{m}}, \ \text{where} \ V_{\underline{m}} = \{v \in V \ | \ d_i.v = m_iv, \ 0 \leqslant i \leqslant n \}$;
	\item $\dim V_{\underline{m}} < \infty \ \forall \ \underline{m} \in \mathbb{Z}^{n+1}$.
\end{enumerate}  
\noindent The collection $P(V) = \{\mu \in D^* \ | \ V_{\mu} \neq (0) \}$, where $V_{\mu} = \{v \in V \ | \ d_iv = \mu(d_i)v, \ 0 \leqslant i \leqslant n \}$ and $D=$ span$\{d_0, \ldots, d_n\}$, is known as the set of all weights of $V$ with respect to $D$. Elements of $V_{\mu}$ are said to be  weight vectors of weight $\mu$.
\brmk\label{R2.2}
If $V$ is an irreducible (and hence indecomposable) Harish-Chandra module over $\tau$, then it is easy to see that there exists $\lambda \in D^{*}$ such that $P(V) \subseteq \{\lambda + \underline{m} \ | \ \underline{m} \in \mathbb{Z}^{n+1} \}$.
\ermk

\section{Classification of Irreducible Cuspidal Modules}\label{Central}

In this section, we introduce the notion of cuspidal modules over $\tau$ and then classify all the irreducible cuspidal modules up to isomorphism.
\bdfn
A Harish-Chandra module $V$ over $\tau$ is said to be cuspidal if there
exists $N_0 \in \mathbb{N}$ such that $\dim V_{\underline{m}} \leqslant N_0 \ \forall \ \underline{m} \in \mathbb{Z}^{n+1}$.
\edfn 
The following lemma is elementary, but we provide a proof for the convenience of the reader, as it will be used repeatedly throughout this paper.
\blmma\label{L2.2}
Let $L$ be any Lie algebra and $V$ be  any non-zero $L$-module. Then $V$ has a non-zero irreducible $L$-subquotient.
\elmma
\begin{proof} If $V$ is a trivial $L$-module, then the result is clear. Else if $L$ does not act trivially on $V$, there exists a non-zero $v_0 \in V$ such that $U(L) v_0 \neq (0)$. Set $V' = U(L)v_0$. If $V'$ is irreducible over $L$, then we are done. Otherwise produce via Zorn's lemma a maximal $L$-submodule $W'$ of $V'$ that does not contain $v_0$. Then $V'/W'$ is an irreducible $L$-subquotient of $V$.
\end{proof}

\subsection{Action of the Central Operators} \label{S2}
In this subsection, unless otherwise stated, $V$ will always stand for an irreducible cuspidal module over $\tau$. Under this assumption, our initial aim in this subsection is to prove the following result.
\bthm\label{T2.3}
$\mathcal{Z}$ acts trivially on $V$.
\ethm
\noindent We need some preparation to prove this theorem. We first introduce the notion of a central operator for a $\tau$-module, which plays a vital role in the proof.

\bdfn
For any $\tau$-module $V$, a linear map $z : V \longrightarrow V$ is called a $\tau$-central operator of degree $\underline{m}$ if it satisfies the following conditions.
\begin{enumerate}
	\item $z$ commutes with the action of $\overline{L}(\mathfrak{g})$,  
	\item $d_i z -z d_i = m_i z \ \forall \ i = 0, \ldots ,n$.
\end{enumerate}
\edfn

\blmma\label{L2.5}
The central element $K_i$ acts trivially on $V$ for each $0 \leqslant i \leqslant n$.
\elmma 
\begin{proof}
Fix any $0 \leqslant i \leqslant n$. Since $V$ is an irreducible Harish-Chandra module, $K_i$ acts on $V$ by a fixed scalar. Recall that $P(V) \subseteq \{ \lambda + \underline{m} \ | \ \underline{m} \in \mathbb{Z}^{n+1} \}$ for some $\lambda \in D^*$. Assume that $V_\lambda \neq 0$. Define $M := \bigoplus_{m \in \mathbb{Z}} V_{\lambda + m e_i}$, which is a non-zero cuspidal module over the classical Virasoro algebra with center $K_i$ and so $K_i$ acts by zero on $M$ by \cite{MP}. This implies that $K_i$ acts trivially on $V$. 
\end{proof}

\blmma \label{L2.6}
Given $0 \leqslant i \leqslant n$ and $\underline{r} \in \mathbb{Z}^{n+1}$,
$t^{\underline{r}}K_i$ is either injective or locally
nilpotent on $V$. 
\elmma

\begin{proof} Assume that there exists some non-zero $w \in V$ such that
$t^{\underline{r}}K_i.w = 0$. Now we already know that $\mathcal{Z}$ is
central in $\overline{L}(\mathfrak{g})$. Also by application of induction
on $p$, it is not too difficult to verify that
\[
(t^{\underline{r}}K_i)^{p+1} \big ((t^{\underline{s_1}}d_{i_1} \ldots
t^{\underline{s_p}}d_{i_p})w \big ) = 0 \qquad
	\text{for all} \ p \in \mathbb{N}, \ \underline{s_1} \ldots,
	\underline{s_p} \in \mathbb{Z}^{n+1} \ \text{and} \ 0 \leqslant
	i_1, \ldots, i_p \leqslant n.
\]
	The required result then follows from the irreducibility of $V$ over 
	$\tau$.
\end{proof} 
\blmma \label{L2.7}
Fix any $0 \leqslant i \leqslant n$ and suppose that there exists $\underline{m} \in \mathbb{Z}^{n+1}$ with $m_j \neq 0$ for some $0 \leqslant j \neq i \leqslant n$ such that $t^{\underline{m}}K_i$ acts  injectively on $V$. Then $t^{\underline{s}}K_i$ acts injectively on $V$ for all $\underline{s} \in \mathbb{Z}^{n+1}$ satisfying $s_k \neq 0$ for some $0 \leqslant k \neq i \leqslant n$.   
\elmma
\begin{proof}
	Assume for contradiction that there exists $\underline{s}
	\in \mathbb{Z}^{n+1}$ with $s_k \neq 0$ for some $0 \leqslant k
	\neq i \leqslant n$ such that $t^{\underline{s}}K_i$ does not act
	injectively on $V$. By Lemma \ref{L2.6},
	$t^{\underline{s}}K_i$ acts locally nilpotently on $V$.

	\noindent \textbf{Claim}.  $t^{\underline{s}}K_i$ acts
	nilpotently on $V$. \\
	Indeed, since $V$ is cuspidal, there exists $N
	\in \mathbb{N}$ such that $\big ((t^{\text{-}\underline{s}}K_i)
	(t^{\underline{s}}K_i) \big)^NV=(0)$. This gives
	\begin{align*}
		0 = t^{\underline{m} + \underline{s}}d_k \bigg ( \big (
		(t^{\text{-}\underline{s}}K_i)^N (t^{\underline{s}}K_i)^N
		\big )v \bigg ) 
		= &\ Ns_k \bigg ((t^{\underline{m} +
		\underline{2s}}K_i)(t^{\text{-}\underline{s}}K_i)^N \big
		((t^{\underline{s}}K_i)^{N-1}v \big ) \\
		&\ - (t^{\underline{m}}K_i)(t^{\text{-}
		\underline{s}}K_i)^{N-1} \big ((t^{\underline{s}}K_i)^Nv
		\big ) \bigg ) \qquad \forall \ v \in V.   
	\end{align*}
	Applying $t^{\underline{s}}K_i$ on both sides, we get
	$(t^{\underline{m}}K_i)(t^{\text{-}
	\underline{s}}K_i)^{N-1}(t^{\underline{s}}K_i)^{N+1}V = (0)$,
	whence by hypothesis, it directly follows that $(t^{\text{-}
	\underline{s}}K_i)^{N-1}(t^{\underline{s}}K_i)^{N+1}V = (0)$.
	Again applying $t^{\underline{m} + \underline{s}}d_k$ another
	$(N-1)$ times and proceeding similarly as above, we finally
	obtain $(t^{\underline{s}}K_i)^{2N}V = (0)$. Hence the claim. 
	
    Taking $N^{\prime} = 2N$, we thereby obtain
	\begin{align*}
		0=t^{\underline{r_1}}d_k \big ((t^{\underline{s}}K_i)^{N^{\prime}}v \big )=N^{\prime}s_k \bigg (t^{\underline{r_1}+ \underline{s}}K_i \big ((t^{\underline{s}}K_i)^{N^{\prime}-1}v \big ) \bigg ) \ \forall \ \underline{r_1}\in \mathbb{Z}^{n+1}, \ v \in V,
	\end{align*}
	from which it clearly follows that
	$\big (t^{\underline{r_1}+ \underline{s}}K_i
	(t^{\underline{s}}K_i)^{N^{\prime}-1} \big )V=(0)$ for all
	$\underline{r_1} \in \mathbb{Z}^{n+1}.$
	This implies:
	\begin{align*}
		t^{\underline{r_2}}d_k \bigg (s_kt^{\underline{r_1}+ \underline{s}}K_i \big ((t^{\underline{s}}K_i)^{N^{\prime}-1}v \big ) \bigg )=0 \ \forall \ \underline{r_1}, \underline{r_2} \ \text{and} \ v \in V,
	\end{align*}
	from which we directly get
	\begin{align*}
		(N^{\prime}-1)s_{k}^{2} \bigg (t^{\underline{r_1}+ \underline{s}}K_it^{\underline{r_2}+ \underline{s}}K_i  \big ((t^{\underline{s}}K_i)^{N^{\prime}-2}v \big ) \bigg )=0 \ \forall \ \underline{r_1},\ \underline{r_2} \in \mathbb{Z}^{n+1}, \ v \in V.
	\end{align*}
	This immediately reveals that
	\begin{align*}
		\big (t^{\underline{r_1}+ \underline{s}}K_it^{\underline{r_2}+ \underline{s}}K_i (t^{\underline{s}}K_i)^{N^{\prime}-2} \big ) V=(0) \ \forall \ \underline{r_1},\ \underline{r_2} \in \mathbb{Z}^{n+1}.
	\end{align*}
	By repeating the above argument another $(N^{\prime}-2)$ times, we obtain
	\begin{align*}
		(t^{\underline{r_1}+ \underline{s}}K_i \ldots t^{\underline{r_{N^{\prime}}}+ \underline{s}}K_i) V=(0) \ \forall \ \underline{r_1}, \ldots, \underline{r_{N^{\prime}}} \in \mathbb{Z}^{n+1}.
	\end{align*}
	Choose any $\underline{s_1}, \ldots,
	\underline{s_{N^{\prime}}} \in \mathbb{Z}^{n+1}$ and set
	$\underline{r_l}= \underline{s_l}- \underline{s} \ \forall \,
	1 \leqslant l \leqslant N^{\prime}$. This and Lemma \ref{L2.5} give 
	\begin{align*}
		(t^{\underline{s_1}}K_i \ldots
		t^{\underline{s_{N^{\prime}}}}K_i) V=(0)
	\end{align*}
	and so there exists a non-zero vector $v_0 \in V$ such that
	$(t^{\underline{r}}K_i)v_0=0 \ \forall \ \underline{r} \in
	\mathbb{Z}^{n+1}$.
	In particular, $t^{\underline{m}}K_i$ does not act injectively on
	$V$, which contradicts our hypothesis and hence proves the
	lemma.
\end{proof}

\blmma \label{L2.8}
Fix any $0 \leqslant i \leqslant n$. Then $t^{\underline{m}}K_i$ cannot
act injectively on $V$ for all $\underline{m} \in \mathbb{Z}^{n+1}$
satisfying $m_j \neq 0$ for some $0 \leqslant j \neq i \leqslant n$.  
\elmma

\begin{proof}
If possible, let $t^{\underline{m}}K_i$ act injectively on $V$ for all
$\underline{m} \in \mathbb{Z}^{n+1}$ satisfying $m_j \neq 0$ for some $0
\leqslant j \neq i \leqslant n$. This implies that $t^{\underline{m}}K_i
: V_{\underline{0}} \longrightarrow V_{\underline{m}}$ is
injective for all $\underline{m} \in \mathbb{Z}^{n+1}$ satisfying
$m_j \neq 0$ for some $0 \leqslant j \neq i \leqslant n$. Thus we
obtain an injective linear operator $T : V_{\underline{0}}
\longrightarrow V_{\underline{m}}$ for all  $\underline{m} \in
\mathbb{Z}^{n+1}$.
Similarly we also have an injective linear operator $T^{\prime} :
V_{\underline{m}} \longrightarrow V_{\underline{0}}$ for all
$\underline{m} \in \mathbb{Z}^{n+1}$. Hence we obtain dim $V_{\underline{0}} =$ dim
$V_{\underline{m}}$ for all $\underline{m} \in \mathbb{Z}^{n+1}$, which
yields that 
	\begin{align}\label{dimension}
	\dim V_{\underline{r}} = \dim V_{\underline{s}} \ \forall \
	\underline{r}, \underline{s} \in \mathbb{Z}^{n+1}
	\end{align}   
	Put $L_i = \{ \underline{m} \in \mathbb{Z}^{n+1} \ | \
	t^{\underline{m}}K_i \ \text{acts injectively on} \ V  \}$ and
	let $S$ be the subgroup generated by $L_i$. Then from our initial
	assumption, we have rank $S=n+1$. Thus there exist
	$\underline{s_0}, \ldots, \underline{s_n} \in \mathbb{Z}^{n+1}$
	such that $m_0\underline{s_0}, \ldots, m_n\underline{s_n}$ forms
	a basis of $S$ for some non-zero integers $m_0, \ldots, m_n$.
	Pick $B \in GL(n+1, \mathbb{Z})$ such that
	$B(\underline{s_k})=e_k$ for all $0 \leqslant k \leqslant n$.
	This gives us $B(m_k\underline{s_k}) = m_ke_k$ for all $0
	\leqslant k \leqslant n$ and so, up to a change of co-ordinates,
	we can assume that there exist injective $\tau$-central operators
	$z_0 \ldots, z_n$ on $V$ of respective degrees $m_0e_0, \ldots,
	m_ne_n$. Again due to (\ref{dimension}), $z_k$ is an invertible
	$\tau$-central operator on $V$ for each $0 \leqslant k \leqslant
	n$. Taking $T_k = z_k^{-1}$, we can check
	that $T_k$ is also a $\tau$-central operator on $V$ of
	degree $-m_ke_k$ for every $0 \leqslant k \leqslant n$.
	Henceforth arguing similarly as in Claim 1, Claim 2 and Claim 3
	of \cite[Theorem 4.5]{E1} and using the fact dim
	$V_{\underline{m}} < \infty$ for all $\underline{m} \in
	\mathbb{Z}^{n+1}$, we can deduce that
	\begin{align*}
		W = \text{span} \{z_kv - v \ | \ v \in V, \ 0 \leqslant k \leqslant n \}
	\end{align*} 
	is a proper $\overline{L}(\mathfrak{g})$-submodule of $V$ with dim $V/W < \infty$. Set $\overline{V} = V/W$.      
	Now note that $ \mathfrak{h} \otimes A \oplus \mathcal{Z}$ is a
	solvable Lie algebra acting on $\overline{V}$. Thus
	by Lie's theorem, there exist a non-zero $\overline{v^{\prime}}
	\in \overline{V}$ and a linear functional $\mu$ such that $(h
	\otimes t^{\underline{r}})\overline{v^{\prime}} = \mu(h,
	\underline{r})\overline{v^{\prime}}$ for all $h \in \mathfrak{h}$
	and $\underline{r} \in \mathbb{Z}^{n+1}$, whence we get
	$\mathcal{Z}\overline{v^{\prime}} = [\mathfrak{h} \otimes A,
	\mathfrak{h} \otimes A]\overline{v^{\prime}} = \overline{0}$. In
	particular, $z_0 \overline{v^{\prime}} = \overline{0}$ and
	so $z_0v^{\prime} \in W$. Moreover $T_0(z_kv - v) =
	z_k(T_0v) - T_0v \in W \ \forall \ v \in V, \ 0 \leqslant k
	\leqslant n$. This finally implies that
	$v^{\prime} = T_0(z_0v^{\prime}) \in T_0(W) \subseteq W,$
	which is a contradiction and hence the lemma is proved.
\end{proof}

\noindent \textbf{Proof of Theorem \ref{T2.3}}.
Lemmas \ref{L2.6}, \ref{L2.7} and \ref{L2.8} reveal that
$t^{\underline{k}}K_i$ acts locally nilpotently on $V$ for all $0
\leqslant i \leqslant n$ and $\underline{k} \in \mathbb{Z}^{n+1}$
with $k_j \neq 0$ for some $0 \leqslant j \neq i \leqslant n$.
Since $V$ is cuspidal, there exists $N \in \mathbb{N}$ such
that $\big ( (t^{\underline{k}}K_i)(t^{\text{-}\underline{k}}K_i)
\big)^NV=0$ for all $\underline{k} \in \mathbb{Z}^{n+1}$
with $k_j \neq 0$ for some $0 \leqslant j \neq i \leqslant n$. Then
for any $\underline{m} \in \mathbb{Z}^{n+1}$,    
\begin{align*}
	0 = t^{\underline{m}}d_j \bigg ( \big ( (t^{\underline{k}}K_i)^N
	(t^{\text{-}\underline{k}}K_i)^N \big )v \bigg ) = &\ N k_j \bigg
	((t^{\underline{m} +
	\underline{k}}K_i)(t^{\underline{k}}K_i)^{N-1} \big
	((t^{\text{-}\underline{k}}K_i)^Nv \big ) \\
	&\ - (t^{\underline{m} - \underline{k}}K_i)(t^{\text{-}
	\underline{k}}K_i)^{N-1} \big ((t^{\underline{k}}K_i)^Nv \big )
	\bigg ) \qquad \forall \ v \in V.
\end{align*}
Applying $t^{\underline{k}}K_i$ on both sides, we thereby obtain 
\begin{align*}
	\big ((t^{\underline{m} - \underline{k}}K_i)(t^{\text{-} \underline{k}}K_i)^{N-1} (t^{\underline{k}}K_i)^{N+1} \big )V = (0)
\end{align*}
for all $\underline{m} \in \mathbb{Z}^{n+1}$ and for each $\underline{k} \in \mathbb{Z}^{n+1}$ with $k_j \neq 0$ for some $0 \leqslant j \neq i \leqslant n$.\\
This implies that for every $\underline{m_1}, \underline{m_2}  \in \mathbb{Z}^{n+1}$, we have   
\begin{align*}
0 = &\ t^{\underline{m_2}}d_j \bigg ((t^{\underline{m_1} -
\underline{k}}K_i)(t^{\text{-} \underline{k}}K_i)^{N-1}
\big((t^{\underline{k}}K_i)^{N+1}v \big ) \bigg ) \\
= &\ k_j \bigg ((N+1)(t^{\underline{m_2} +
\underline{k}}K_i)(t^{\underline{m_1} - \underline{k}}K_i)(t^{\text{-}
\underline{k}}K_i)^{N-1} \big ((t^{\underline{k}}K_i)^Nv \big ) \\
&\ - (N-1)(t^{\underline{m_2} - \underline{k}}K_i)(t^{\underline{m_1} -
\underline{k}}K_i)(t^{\text{-} \underline{k}}K_i)^{N-2} \big
((t^{\underline{k}}K_i)^{N+1}v \big ) \bigg) \qquad \forall \ v \in V.
\end{align*}   
Again applying $t^{\underline{k}}K_i$ on both sides, we get
$(t^{\underline{m_1} - \underline{k}}K_i)(t^{\underline{m_2} -
\underline{k}}K_i)(t^{\text{-} \underline{k}}K_i)^{N-2} \big
(t^{\underline{k}}K_i)^{N+2}V = (0)$
for all $\underline{m_1}, \underline{m_2}  \in \mathbb{Z}^{n+1}$ and each $\underline{k} \in \mathbb{Z}^{n+1}$ with $k_j \neq 0$ for some $0 \leqslant j \neq i \leqslant n$.\\
Repeating the above argument another $(N-2)$ times gives
\begin{align*}
	(t^{\underline{m_1} - \underline{k}}K_i)(t^{\underline{m_2} - \underline{k}}K_i) \ldots (t^{\underline{m_N} - \underline{k}}K_i) (t^{\underline{k}}K_i)^{2N}V = (0)
\end{align*}

\noindent for all $\underline{m_1}, \underline{m_2}, \ldots,
\underline{m_N}  \in \mathbb{Z}^{n+1}$ and for every $\underline{k} \in
\mathbb{Z}^{n+1}$ satisfying $k_j \neq 0$ for some $0 \leqslant j \neq i
\leqslant n$. Consequently, we have 
\begin{align*}
0 = &\ (t^{\underline{m_{N+1}} - 2 \underline{k}}d_j) \bigg
((t^{\underline{m_1} - \underline{k}}K_i)(t^{\underline{m_2} -
\underline{k}}K_i) \ldots (t^{\underline{m_N} - \underline{k}}K_i) \big
((t^{\underline{k}}K_i)^{2N}v \big ) \bigg ) \\
= &\ Nk_j \bigg ((t^{\underline{m_1} -
\underline{k}}K_i)(t^{\underline{m_2} - \underline{k}}K_i) \ldots
(t^{\underline{m_{N+1}} - \underline{k}}K_i) \big
((t^{\underline{k}}K_i)^{2N-1}v \big ) \\
&\ + (t^{\underline{m_1} - \underline{k}}K_i)(t^{\underline{m_2} -
\underline{k}}K_i) \ldots (t^{\underline{m_{N+1}} - \underline{k}}K_i)
\big ((t^{\underline{k}}K_i)^{2N-2}v \big ) \bigg ) \\
&\ + (t^{\underline{k}}K_i)^{2N-1}(t^{\underline{m_1} -
\underline{k}}K_i) \ldots (t^{\underline{m_N} - \underline{k}}K_i) \bigg
( k_j \big (t^{\underline{m_{N+1}} - \underline{k}}K_i)v \big ) +
(t^{\underline{k}}K_i) \big ((t^{\underline{m_{N+1}} - 2
\underline{k}}d_j)v \big ) \bigg ) \ \forall \ v \in V.
\end{align*}
Applying $t^{\underline{k}}K_i$ on both sides of the above equation, we obtain 
\begin{align*}
	(t^{\underline{m_1} - \underline{k}}K_i)(t^{\underline{m_2} - \underline{k}}K_i) \ldots (t^{\underline{m_{N+1}} - \underline{k}}K_i) (t^{\underline{k}}K_i)^{2N-1}V = (0)
\end{align*}
for all $\underline{m_1}, \underline{m_2}, \ldots, \underline{m_{N+1}}  \in \mathbb{Z}^{n+1}$ and for every $\underline{k} \in \mathbb{Z}^{n+1}$ satisfying $k_j \neq 0$ for some $0 \leqslant j \neq i \leqslant n$.\\  
Repeating the above argument another $(2N-1)$ times, we get
\begin{align*}
	(t^{\underline{m_1} - \underline{k}}K_i)(t^{\underline{m_2} - \underline{k}}K_i) \ldots (t^{\underline{m_{3N}} - \underline{k}}K_i) V = (0)
\end{align*}
for all $\underline{m_1}, \underline{m_2}, \ldots, \underline{m_{3N}}  \in \mathbb{Z}^{n+1}$ and for every $\underline{k} \in \mathbb{Z}^{n+1}$ satisfying $k_j \neq 0$ for some $0 \leqslant j \neq i \leqslant n$.\\
Let $\underline{k_1}, \ldots, \underline{k_{3N}} \in \mathbb{Z}^{n+1}$ be arbitrary. Choose any $\underline{k} \in \mathbb{Z}^{n+1}$ such that $k_j \neq 0$ for some $0 \leqslant j \neq i \leqslant n$  and thereby set $\underline{m_l}= \underline{k_l} + \underline{k} \ \forall \ l=1, \ldots, 3N$.  Along with Lemma \ref{L2.5}, this ultimately gives us
\begin{align*}
	(t^{\underline{k_1}}K_i \ldots t^{\underline{k_{3N}}}K_i) V = (0) \ \forall \ i=0, \ldots, n.
\end{align*}
Subsequently by proceeding exactly as in \cite[Lemma 4.4]{PR}, we can  find a non-zero vector $v_0 \in V$ such that $(t^{\underline{k}}K_i)v_0=0 \ \forall \ \underline{k} \in \mathbb{Z}^{n+1}, \ 0 \leqslant i \leqslant n$. Now it is trivial to check that $W= \{v \in V \ | \ \mathcal{Z}v = (0) \}$ 
is a $\tau$-submodule of $V$. Therefore we are done by the irreducibility of $V$.
\qed

Our next goal is to establish that $\mathcal{Z}$ acts trivially on every
irreducible cuspidal module for the higher-dimensional Virasoro algebra
$\mathcal{V}ir$. Since $\mathcal{V}ir$ is a subalgebra of $\tau$,
the techniques used in proving Theorem \ref{T2.3} and Lemmas
\ref{L2.5}, \ref{L2.6} and \ref{L2.7} also apply to irreducible
cuspidal modules over $\mathcal{V}ir$. However, this is not the case for
the proof of Lemma \ref{L2.8}, as we have used $\mathcal{Z} =
[\mathfrak{h} \otimes A , \mathfrak{h} \otimes A]$. So to prove
Theorem \ref{T2.3} for $\mathcal{V}ir$-modules, it suffices to
prove the analogue of Lemma \ref{L2.8} for them.

To this end, we first define the notion of central operators for
$(\mathcal{V}ir)_{0}$-modules, where  $(\mathcal{V}ir)_{0} = \text{span}
\{t^{\underline{r}}d_i, t^{\underline{s}}K_j \ | \ \underline{r},
\underline{s} \in \{\underline{0}\} \times \mathbb{Z}^n, 0 \leqslant i,j \leqslant n \}$ is a
subalgebra of $\mathcal{V}ir$.

\bdfn
For any $(\mathcal{V}ir)_{0}$-module $V$, a linear map $z : V \longrightarrow V$ is said to be a $(\mathcal{V}ir)_{0}$-central operator of degree $\underline{m}$ if it satisfies the following conditions.
\begin{enumerate}
	\item $z$ commutes with the action of $(\mathcal{V}ir)_{0}^{\prime}$, where \\ $(\mathcal{V}ir)_{0}^{\prime} = \text{span} \{ t^{\underline{p}}K_0, t^{\underline{q}}d_0, t^{\underline{r}}K_j \ | \ \underline{p}, \underline{q}, \underline{r} \in \{\underline{0}\} \times \mathbb{Z}^n, 1 \leqslant j \leqslant n \}$ and  
	\item $d_i z -z d_i = m_i z \ \forall \ i = 1, \ldots ,n$.
\end{enumerate}
\edfn
\brmk
$t^{\underline{m}}K_j$ is a $(\mathcal{V}ir)_{0}$-central operator for any $\underline{m} \in \{\underline{0}\} \times \mathbb{Z}^n$ and $1 \leqslant j \leqslant n$.
\ermk
\blmma \label{L2.12}
Let $V$ be an irreducible cuspidal module over $\mathcal{V}ir$.
Then for any $0 \leqslant i \leqslant n$, $t^{\underline{m}}K_i$ cannot act injectively on $V$ for all $\underline{m} \in \mathbb{Z}^{n+1}$ satisfying $m_j \neq 0$ for some $0 \leqslant j \neq i \leqslant n$.  
\elmma
\begin{proof}
For the sake of contradiction, suppose that $t^{\underline{m}}K_i$ acts injectively on $V$ for all
$\underline{m} \in \mathbb{Z}^{n+1}$ satisfying $m_j \neq 0$ for some $0
\leqslant j \neq i \leqslant n$. Without loss of generality, let us take $i=0$.
Thus $t^{\underline{m}}K_0$ acts injectively on $V$ for every
$\underline{0} \neq \underline{m} = (0,m_1, \ldots, m_n) \in \mathbb{Z}^{n+1}$.
Now use Lemma \ref{L2.2} to obtain an irreducible
$(\mathcal{V}ir)_{0}$-subquotient of $V$, say $\overline{V}$. Then
$\overline{V}$ is a cuspidal module over $(\mathcal{V}ir)_{0}$ with
respect to $D^{\prime} = \text{span} \{d_1, \ldots, d_n \}$, since $d_0$
is central in $(\mathcal{V}ir)_{0}$.

\noindent \textbf{Case 1.} $n =1$. \\
By Lemma \ref{L2.5}, $\mathcal{H} := \text{HVir}/ \text{span} \{C_D, C_{DI}, C_I\} \cong
\text{span} \{t_1^{r}d_1, t_1^{s}K_0 \ | \ r, s \in \mathbb{Z} \}$ is a
subalgebra of $(\mathcal{V}ir)_0$ (see Subsection \ref{Twisted}). So
$\overline{V}$ is a cuspidal $\mathcal{H}$-module that contains the rank
$1$ Witt algebra $\mathcal{W}_1 = \text{span} \{t_1^{r}d_1 \ | \ r \in
\mathbb{Z} \}$ as a subalgebra. Thus by Remark \ref{R1.1} and
\cite[Lemma 3.3]{Ma}, $\overline{V}$ has finite length and so
contains an irreducible $\mathcal{H}$-submodule, say
$\overline{W} = W_1/W_2$, where $W_1$ and $W_2$ are subspaces of $V$.
But since $K_0$ acts trivially on $\overline{W}$ (by Lemma
\ref{L2.5}), we can apply \cite[Theorem 4.4]{LZ2} to
conclude that $t_1^{s}K_0$ acts trivially on $\overline{W}$ for all
$s \in \mathbb{Z}$. This contradicts our initial assumption.

\noindent \textbf{Case 2.} $n \geqslant 2$.
We split this proof into the following three steps.\smallskip

\noindent \textbf{Step 1.} $t^{\underline{m}}K_i$ acts locally
nilpotently on $\overline{V}$ for all $\underline{m} \in \{\underline{0}\} \times \mathbb{Z}^n$
and each $1 \leqslant i \leqslant n$.

\noindent \textbf{Step 2.} $\mathcal{Z}^{\prime} = \text{span}
\{t^{\underline{r}}K_j \ | \ \underline{r} \in \{\underline{0}\} \times \mathbb{Z}^n, 1 \leqslant
j \leqslant n \}$ acts trivially on $\overline{V}$.

\noindent \textbf{Step 3.} For each $\underline{m} \in \{\underline{0}\} \times \mathbb{Z}^n$,
$t^{\underline{m}}K_0$ acts trivially on a non-zero vector space
$\overline{M} = M_1/M_2$, where $M_1$ and $M_2$ are subspaces of
$V$.\medskip
	
	$Step \ 1.$ If the claim of Step 1 is false, then as $n
	\geqslant 2$, we can proceed similarly as in Lemma \ref{L2.7} to
	show that there exists $1 \leqslant i \leqslant n$ such that 
	$t^{\underline{m}}K_i$ is injective on $\overline{V}$ for
	all $\underline{0} \neq \underline{m} \in \{\underline{0}\} \times \mathbb{Z}^{n}$. This
	implies $\dim \overline{V}_{\underline{r}} = \dim
	\overline{V}_{\underline{s}} \ \forall \ \underline{r},
	\underline{s} \in \{\underline{0}\} \times \mathbb{Z}^n$. Then our assumption
	yields that $t^{\underline{m}}K_i$ is invertible for each
	$\underline{0} \neq \underline{m} \in \{\underline{0}\} \times \mathbb{Z}^n$.
	So as in Lemma \ref{L2.8}, there exist invertible
	$(\mathcal{V}ir)_{0}$-central operators $z_1, \ldots, z_n$ on
	$\overline{V}$ with degrees $l_1e_1, \ldots, l_ne_n$ for non-zero $l_1, \ldots, l_n \in \mathbb{Z}$. By providing a
	similar argument as in Claims 1, 2, and 3 of \cite[Theorem
	4.5]{E1} and using $\dim \overline{V}_{\underline{m}}
	< \infty\ \forall \ \underline{m} \in \{\underline{0}\} \times \mathbb{Z}^{n}$, we infer
	\begin{align*} \overline{W} = \text{span} \{z_k \overline{v}
	- \overline{v} \ | \ \overline{v} \in \overline{V}, \ 1 \leqslant
	k \leqslant n \}
	\end{align*} 
	is a proper $(\mathcal{V}ir)_{0}^{\prime}$-submodule of
	$\overline{V}$ with dim $\overline{V}/ \overline{W} < \infty$.
	Again note that $(\mathcal{V}ir)_{0}^{\prime}$ is solvable
	and 
	\begin{align*}
		\bigg[\bigoplus_{\underline{r} \in \{\underline{0}\} \times \mathbb{Z}^{n}}\mathbb{C}t^{\underline{r}}d_0, \bigoplus_{\underline{s} \in \{\underline{0}\} \times \mathbb{Z}^{n}}\mathbb{C}t^{\underline{s}}K_0 \bigg] = \mathcal{Z}^{\prime},    
	\end{align*} 
	which ultimately leads to a contradiction, by appealing
	essentially to the same argument that we have used in Lemma
	\ref{L2.8}. Hence the first step follows.\smallskip
	
	$Step \ 2.$ As $n \geqslant 2$, we can use Step 1 and then proceed similarly as in Theorem \ref{T2.3} to obtain $\overline{0} \neq \overline{v_0} \in \overline{V}$ such that $t^{\underline{m}}K_j \overline{v_0} = \overline{0} \ \forall \ \underline{m} \in \{\underline{0}\} \times \mathbb{Z}^{n}, \ 1 \leqslant j \leqslant n$. This  shows that $\overline{T}= \{v \in \overline{V} \ | \ \mathcal{Z}^{\prime}\overline{v} = (\overline{0}) \}$ is a non-zero $(\mathcal{V}ir)_{0}$-module. The second step now follows from the irreducibility of
	$\overline{V}$.\smallskip
	
	$Step \ 3.$ By Step 2, $\overline{V}$ is cuspidal over $(\mathcal{V}ir)_{0}^{\prime \prime} = \text{span} \{t^{\underline{r}}K_0, t^{\underline{s}}d_j \ | \ \underline{r}, \underline{s} \in \{\underline{0}\} \times \mathbb{Z}^{n}, 1 \leqslant j \leqslant n \}$. Furthermore the rank $n$ Witt algebra $\mathcal{W}_n = \text{span} \{t^{\underline{r}}d_j \ | \ \underline{r} \in \{\underline{0}\} \times \mathbb{Z}^n, 1 \leqslant j \leqslant n \}$ also sits inside $(\mathcal{V}ir)_{0}^{\prime \prime}$ as a subalgebra. Consequently \cite[Lemma 3.3]{Ma} and Remark \ref{R1.1} reveals that $\overline{V}$ has a finite composition series, which gives rise to a non-zero irreducible cuspidal $(\mathcal{V}ir)_{0}^{\prime \prime}$-submodule of $\overline{V}$, say $\overline{M} = M_1/M_2$, where $M_1$ and $M_2$ are subspaces of $V$. Finally as $n \geqslant 2$ and $K_0$ acts trivially on $\overline{M}$ (by Lemma \ref{L2.5}), we can now directly apply \cite[Theorem 3.3]{GLZ} to obtain the third step.
	
This is a contradiction to our initial assumption, which thereby gives us the desired result.
\end{proof}

\bthm \label{T2.13}
$\mathcal{Z}$ acts trivially on every irreducible cuspidal module
over $\mathcal{V}ir$ (or over $\tau$). 
\ethm 
\begin{proof}
If $V$ is an irreducible cuspidal $\mathcal{V}ir$-module, the
result follows by combining Lemmas \ref{L2.7} and \ref{L2.12} and
Theorem \ref{T2.3}. For an irreducible cuspidal $\tau$-module $V$,
this is precisely Theorem \ref{T2.3}. 
\end{proof}

\brmk\
The above theorem reduces our problem of classifying irreducible cuspidal modules over $\tau$ to classifying all those irreducible cuspidal modules for $\widehat{\tau}$, where $\widehat{\tau} = \mathcal{W}_{n+1} \ltimes L(\mathfrak{g})$.
\ermk

\subsection{Jet Modules} 
Set $A \widehat{\tau} = \widehat{\tau} \oplus A$. Then $A \widehat{\tau}$ forms a Lie algebra by simply extending the Lie algebra structure already prevalent on $\widehat{\tau}$ in the following manner (see \cite{Y,E2}). 
\begin{align*}
[t^{\underline{m}}d_i, t^{\underline{k}}] = k_it^{\underline{k} + \underline{m}}d_i, \,\,\,\
[t^{\underline{k}}, x \otimes t^{\underline{m}}] = 0
\end{align*} 
for all $x \in \mathfrak{g}, \ \underline{k}, \underline{m} \in \mathbb{Z}^{n+1}, \ 0 \leqslant i \leqslant n$.\\
We now recall the following actions, which will be used later in this section.\\
(A1) $A$ acts on $\mathcal{W}_{n+1}$ via 
$t^{\underline{k}} (t^{\underline{m}}d_i) = t^{\underline{k} + \underline{m}}d_i \ \forall \ \underline{k}, \underline{m} \in \mathbb{Z}^{n+1}, \ 0 \leqslant i \leqslant n$. \\
(A2) $A$ acts on $L(\mathfrak{g})$ via $t^{\underline{k}}(x \otimes t^{\underline{m}}) = x \otimes t^{\underline{k} + \underline{m}} \ \forall \ x \in \mathfrak{g}, \ \underline{k}, \underline{m} \in \mathbb{Z}^{n+1}$. \\ 
(A3) $\mathcal{W}_{n+1}$ acts on $A$ via $t^{\underline{m}}d_i( t^{\underline{k}}) = k_it^{\underline{k} + \underline{m}}d_i \ \forall \ \underline{k}, \underline{m} \in \mathbb{Z}^{n+1}, \ 0 \leqslant i \leqslant n$.
\bdfn
An $A\widehat{\tau}$-module $V$ is called a Jet module for $\widehat{\tau}$ if 
\begin{enumerate}
\item $V$ is a Harish-Chandra module with respect to $D = \text{span} \{d_0, \ldots, d_n \}$.
\item $A$ acts associatively on $V$, in the sense that\\
 $1v = v$ and $t^{\underline{r}}(t^{\underline{s}}v) = t^{\underline{r} + \underline{s}}v \ \forall \ \underline{r}, \underline{s} \in \mathbb{Z}^{n+1}, \ v \in V$.
\item $t^{\underline{r}} \big ((x \otimes t^{\underline{s}})v \big ) = (x \otimes t^{\underline{s}})(t^{\underline{r}}v) \  \forall \ x \in \mathfrak{g}, \ \underline{r}, \underline{s} \in \mathbb{Z}^{n+1}, \ v \in V$.     
\end{enumerate}
\edfn

\bthm\cite[Theorem 5.1(iii)]{Y} \label{Jets}
Let $V$ be an irreducible Jet module for $\widehat{\tau}$. Then there exist a finite-dimensional irreducible $\mathfrak{g}$-module $V_1$, a finite-dimensional irreducible $\mathfrak{gl}_{n+1}$-module $V_2$ and some $\underline{\alpha} = (\alpha_0, \ldots, \alpha_n) \in \mathbb{C}^{n+1}$  such that 
\begin{align*}
	V \cong V_1 \otimes V_2 \otimes A,
\end{align*}	
where the action of $\widehat{\tau}$ on $V$ is given by 
\begin{align*}
(x \otimes t^{\underline{m}}) (v_1 \otimes v_2 \otimes t^{\underline{r}})
= &\ (xv_1) \otimes v_2 \otimes t^{\underline{r} + \underline{m}}, \\ 
t^{\underline{m}}d_i (v_1 \otimes v_2 \otimes t^{\underline{r}}) = &\
(\alpha_i + r_i)(v_1 \otimes v_2 \otimes t^{\underline{r} +
\underline{m}}) + \sum_{j=0}^{n}m_j \big (v_1 \otimes (E_{j,i}v_2)
\otimes t^{\underline{r} + \underline{m}} \big ) \\
&\ \forall \ x \in \mathfrak{g}, \ v_1 \in V_1, \ v_2 \in V_2, \
\underline{m},\underline{r} \in \mathbb{Z}^{n+1}, \ 0 \leqslant i
\leqslant n.    
\end{align*} 
Here $E_{j,i}$ denotes the matrix of order $(n+1)$ having $1$ at the $(j,i)$-th entry and $0$ elsewhere.
\ethm
\brmk \label{Quadruplet}If $\mathfrak{g}_1$ is a finite-dimensional simple Lie algebra with a Cartan subalgebra $\mathfrak{h}_1$ and $\{\alpha_1^{\vee}, \cdots ,\alpha_l^{\vee}\}$ is a collection of simple co-roots of $\mathfrak{g}_1$ with respect to $\mathfrak{h}_1$, then we shall denote the set of all dominant integral weights of $\mathfrak{g}_1$ relative to $\mathfrak{h}_1$  by $P_{\mathfrak{g}_1}^+ = \{ \lambda\ \in \mathfrak{h}_1^*\ | \ \lambda(\alpha_i^{\vee}) \in \mathbb{Z}_{+}\ \forall \ i=1, \cdots ,l \}$.   
Now it is a well-known fact that $V_1 \cong V(\lambda_1)$ as  $\mathfrak{g}$-modules for a unique $\lambda_1 \in P_{\mathfrak{g}}^{+}$ and $V_2 \cong V(c, \lambda_2)$ as $\mathfrak{gl}_{n+1}$-modules for a unique pair $(c, \lambda_2) \in \mathbb{C} \times P_{\mathfrak{sl}_{n+1}}^{+}$. So the irreducible module in Theorem \ref{Jets} is completely determined by the quadruplet $(c, \lambda_1, \lambda_2, \underline{\alpha}) \in \mathbb{C} \times P_{\mathfrak{g}}^{+} \times P_{\mathfrak{sl}_{n+1}}^{+} \times \mathbb{C}^{n+1}$ and thus we shall denote this module by $V(c, \lambda_1, \lambda_2, \underline{\alpha})$. Finally it is worth pointing out that 
although $V(c, \lambda_1, \lambda_2, \underline{\alpha})$ is always irreducible as an $A\widehat{\tau}$-module, it need not be irreducible as a $\widehat{\tau}$-module without the $A$-action. 
\ermk

\subsection{Modules of Tensor Fields}

Consider the finite-dimensional simple Lie algebra $\mathfrak{sl}_{n+1}$
and let $H = \text{span} \{\alpha_i^{\vee} := E_{i-1,i-1} -
E_{i,i} \ | \ 1 \leqslant i \leqslant n\}$ be its Cartan subalgebra.
Let $H^*$ denote the dual space of $H$ and $\omega_1, \ldots, \omega_n$
be the $fundamental$ $weights$ of $\mathfrak{sl}_{n+1}$, which is defined
as $\omega_i(\alpha_j^{\vee}) = \delta_{ij}$ for all $1 \leqslant i,j
\leqslant n$. For notational convenience,  we shall also take $\omega_0 =
\omega_{n+1} =0$.

The modules of the form $V(c, 0 , \lambda_2, \underline{\alpha})$ that are $irreducible$ over $\mathcal{W}_{n+1}$ are modules of $tensor$ $fields$ on an $(n+1)$-$dimensional$ $torus$. These modules have their roots in geometry and they are also referred to in the literature as $Larsson$ modules or $Shen$ modules \cite{L,GS}. We now recall the description of the irreducible $\mathcal{W}_{n+1}$-submodules of $V(c, 0, \lambda_2, \underline{\alpha})$ from \cite{BF,E4}, which is essential for our classification of irreducible cuspidal modules.

Let $W \cong \mathbb{C}^{n+1}$ be the standard representation of
$\mathfrak{gl}_{n+1}$. One checks that its $k$-fold exterior power
$\wedge^k W \cong V(\omega_k,k)$ as irreducible
$\mathfrak{gl}_{n+1}$-modules for each $0 \leqslant k \leqslant n+1$.
The corresponding $\mathcal{W}_{n+1}$-modules $V(k , 0, \omega_k,
\underline{\alpha})$ consist of differential $k$-forms and they
form the de Rham complex
\begin{align*}
V(0 , 0, \omega_0, \underline{\alpha}) \xrightarrow{\pi_0} V(1 , 0, \omega_1, \underline{\alpha}) \xrightarrow{\pi_1} \ldots \ldots \xrightarrow{\pi_n} V(n+1 , 0, \omega_{n+1}, \underline{\alpha}) 
\end{align*}
The above homomorphisms of the de Rham complex is a $\mathcal{W}_{n+1}$-module  homomorphism and so the kernels and images of $\pi_i$'s naturally give rise to $\mathcal{W}_{n+1}$-submodules of $V(k , 0, \omega_k, \underline{\alpha})$.  

\bppsn \label{Fundamental}
\
\begin{enumerate}
\item $V(c, \lambda_1, \lambda_2, \underline{\alpha})$ is an irreducible module over $\widehat{\tau}$ if either $(\lambda_2,c, \underline{\alpha}) \notin \{0\} \times \{0, n+1 \} \times \mathbb{Z}^{n+1}$ or $\lambda_1\neq 0$ or $(\lambda_2,c) \neq (\omega_k,k)$ for any fundamental weight $\omega_k$ of $\mathfrak{sl}_{n+1}$, with $1 \leqslant k \leqslant n$.

\item For any $0 \leqslant k \leqslant n+1$, $V(k , 0, \omega_k, \underline{\alpha})$ has a unique irreducible quotient over $\mathcal{W}_{n+1}$ and this irreducible quotient is given by $\pi_k \big(V(k , 0,  \omega_k, \underline{\alpha}) \big)$ when $0 \leqslant k \leqslant n$. Moreover if $\underline{\alpha} \in \mathbb{Z}^{n+1}$, then $V(n+1 , 0, \omega_{n+1}, \underline{\alpha})$ has a trivial $1$-dimensional module as an irreducible quotient.       
\end{enumerate}

\eppsn

\begin{proof}
The second part follows from \cite{E4} (also see \cite{GZ}).
To show the first part, first suppose
$\lambda_1 \neq 0$.
By our construction, $V(\lambda_1)$ is an irreducible module for $\mathfrak{g}$,
which implies that $V(\lambda_1) \otimes A$ is an irreducible module over
$L(\mathfrak{g}) \oplus D$, where $D$ = span$\{d_0, \ldots, d_n \}$. Now
since $\lambda_1 \neq 0$, the required result can be deduced by an
application of \cite[Proposition 2.8]{E2}, which directly conveys that
$V(c,\lambda_2) \otimes A$ is irreducible over $\mathcal{W}_{n+1} \rtimes \big
(\mathfrak{h} \otimes  A \big )$, where $\mathfrak{h}$ is a Cartan
subalgebra of $\mathfrak{g}$.

If instead $\lambda_1 = 0$, then $V(c,
\lambda_1, \lambda_2, \underline{\alpha})$ is an irreducible
$\widehat{\tau}$-module if and only if it is irreducible over
$\mathcal{W}_{n+1}$. The proposition now directly follows from  \cite[Proposition 4.1]{E5} and \cite[Theorem 1.9]{E4} (also see
\cite{GZ}).
\end{proof}

\subsection{Cuspidal Cover} 
Throughout this subsection, unless otherwise stated, $M$ will always stand for an indecomposable cuspidal module over $\widehat{\tau}$. Put
\[\ {\mathcal{L}}= 
\begin{cases}
	\mathcal{W}_{n+1}\ , &\text{if} \,\,\,\,L(\mathfrak{g}) \ \text{acts  trivially on} \ M \\
	L(\mathfrak{g})  \,\,\,\ ,           & \text{otherwise}.
\end{cases}
\]

\bdfn
An $A$-cover of a $\widehat{\tau}$-module $M$ is the subspace
$\widehat{M}$ of $\text{Hom}_\mathbb{C}(A,M)$ that is spanned by
the elements of the form $\mu_{x,u} \ \forall \ x \in \mathcal{L}, \ u
\in M$, with $\mu_{x,u}(a) = (ax)u \ \forall \ a \in A$. 
\edfn

\bppsn \label{Cover}
\
\begin{enumerate}
\item $\widehat{M}$ is an $A\widehat{\tau}$-module under the following actions.\\
$(i) \ (x^{\prime} \otimes t^{\underline{k}})\mu_{x,u} = \mu_{[x^{\prime} \otimes t^{\underline{k}}, x], u} + \mu_{x, (x^{\prime} \otimes t^{\underline{k}})u}, \\
(ii) \ (t^{\underline{k}}d_i)\mu_{x,u} = \mu_{[t^{\underline{k}}d_i, x], u} + \mu_{x, (t^{\underline{k}}d_i)u}, \\ 
(iii) \ (t^{\underline{k}})\mu_{x,u} = \mu_{(t^{\underline{k}})x,u} \ \forall \ x^{\prime} \in \mathfrak{g}, \ \underline{k} \in \mathbb{Z}^{n+1}, \ 0 \leqslant i \leqslant n$.
\item $\widehat{M}$ is a weight module with respect to $D$.
\item There exists a homomorphism of $\widehat{\tau}$-modules $\pi : \widehat{M} \longrightarrow M$ satisfying $\pi(\widehat{M}) = \mathcal{L}M$.  
\end{enumerate}
\eppsn

\begin{proof} (1) This is easily checked using (A1), (A2), (A3) and the bracket operations on $A\widehat{\tau}$.\\
(2) By hypothesis, $M = \bigoplus_{\underline{k} \in \mathbb{Z}^{n+1}} M_{\underline{k}}$. Then for any $u \in M_{\underline{k}}$, we have
\[
d_i \big(\mu_{x^{\prime} \otimes t^{\underline{m}},u} \big) = (k_i + m_i) \big (\mu_{x^{\prime} \otimes t^{\underline{m}},u} \big ), \
d_i \big(\mu_{t^{\underline{m}}d_j, u} \big) = (k_i + m_i) \big
(\mu_{t^{\underline{m}}d_j,u} \big )
\]
for all $x^{\prime} \in \mathfrak{g}, \ \underline{m} \in
\mathbb{Z}^{n+1}, \ 0 \leqslant i,j \leqslant n$.
Now we are done since $\widehat{M}$ is spanned by its weight vectors.

\noindent (3) Define $\pi : \widehat{M} \longrightarrow M$ via
$\pi \big (\mu_{x,u} \big) = \mu_{x,u}(1)$.
Then $\pi$ is a map of $\widehat{\tau}$-modules with
$\pi(\widehat{M}) = \mathcal{L}M$.   
\end{proof}

\blmma \cite[Proposition 5.6]{YV} \label{Differentiators}
Let $N$ be a cuspidal module over $W(\underline{\gamma})$ with respect to $D(\underline{\gamma})$ (see Subsection \ref{Solenoidal}), whose weights are supported on a single coset. Then there exists some $m \in \mathbb{N}$, which depend only on the bound of the dimensions of the weight spaces of $N$, such that 
$\Omega_{\underline{k}, \underline{p}}^{(m, \underline{h})} = \sum_{i=0}^{m}(-1)^i \binom{m}{i} \big (t^{\underline{k} - i \underline{h}}D(\underline{\gamma}) \big) \big (t^{\underline{p} + i \underline{h}}D(\underline{\gamma}) \big)$ in $U(W(\underline{\gamma}))$  annihilates $N$ for all $\underline{k}, \underline{h} \in \mathbb{Z}^{n+1}$ and $\underline{p} \in \mathbb{Z}\underline{h}$. 
\elmma

\brmk \label{Weights}
For a fixed $\lambda \in D^*$ and any $\underline{k} \in \mathbb{Z}^{n+1}$ with $\lambda + \underline{k} \in P(M)$, a direct verification yields that $M_{\lambda + \underline{k}} = M_ {\underline{\alpha}}$, where $\underline{\alpha} = (\underline{\gamma}, \underline{\beta} + \underline{k}), \ \underline{\beta} = (\lambda(d_0), \ldots, \lambda(d_n)) \in \mathbb{C}^{n+1}$ and $(\cdot, \cdot)$ denotes the standard inner product on $\mathbb{C}^{n+1}$. This implies that $M = \bigoplus_{\underline{z} \in \Gamma} M_{\underline{\omega_0} + \underline{z}}$ for a fixed $\underline{\omega_0} = (\underline{\gamma}, \underline{\beta}) \in \mathbb{C}$ and $\Gamma = \{(\underline{\gamma}, \underline{m}) \ | \ \underline{m} \in \mathbb{Z}^{n+1} \} \subseteq \mathbb{C}$, with $M_ {\underline{\omega_0} + \underline{z}} = \{v \in M \ | \ D(\underline{\gamma})v = (\underline{\omega_0} + \underline{z})v \}$.
\ermk

\blmma \label{Cuspidal Cover} 
\
\begin{enumerate}
\item For all $ x \in \mathfrak{g}, \ \underline{j} \in \mathbb{Z}^{n+1}, \ \underline{h} \in \mathbb{Z}^{n+1} \setminus \{\underline{0}\}$ and $\underline{p} \in \mathbb{Z}\underline{h}$, there exists $m \in \mathbb{N}$ such that 
$\mathcal{T}_{\underline{j}, \underline{p}}^{(m, \underline{h})}(x) = \sum_{i=0}^{m}(-1)^i \binom{m}{i} \big (x \otimes t^{\underline{j} - i \underline{h}} \big) \big (t^{\underline{p} + i \underline{h}}D(\underline{\gamma}) \big)$ in $U(\widehat{\tau})$
annihilates $M$.
\item $\widehat{M}$ is a cuspidal module.   
\end{enumerate} 
\elmma

\begin{proof}
(1) By Remark \ref{Weights}, $M$ remains cuspidal when restricted to the solenoidal Lie algebra $W(\underline{\gamma})$ (with respect to $D(\underline{\gamma})$) and all its weights are also supported on a single coset. As a result, we can apply Lemma \ref{Differentiators} to obtain some $r \in \mathbb{N}$ such that for all $x \in \mathfrak{g}, \ \underline{j}, \underline{k}, \ \underline{h} = (h_0, \ldots, h_n) \in \mathbb{Z}^{n+1}$ and $\underline{p} \in \mathbb{Z}\underline{h}$, we have $\big [x \otimes t^{\underline{j}}, \Omega_{\underline{k}, \underline{p}}^{(r, \underline{h})} \big ]M = (0)$. This gives
\begin{align*}
(0) = &\ \bigg (\big [x \otimes t^{\underline{j} + \underline{h}},
\Omega_{\underline{k}, \underline{p} - \underline{h}}^{(r,
\underline{h})} \big] - \big [x \otimes t^{\underline{j}},
\Omega_{\underline{k} + \underline{h}, \underline{p} -
\underline{h}}^{(r, \underline{h})} \big] + 2 \big [x \otimes
t^{\underline{j} - \underline{h}}, \Omega_{\underline{k} + \underline{h},
\underline{p}}^{(r, \underline{h})} \big] \\
&\ - 2 \big [x \otimes t^{\underline{j}}, \Omega_{\underline{k},
\underline{p}}^{(r, \underline{h})} \big] +  \big [x \otimes
t^{\underline{j} - \underline{h}}, \Omega_{\underline{k}, \underline{p} +
\underline{h}}^{(r, \underline{h})} \big] -  \big [x \otimes
t^{\underline{j} - 2\underline{h}}, \Omega_{\underline{k} +
\underline{h}, \underline{p} + \underline{h}}^{(r, \underline{h})} \big]
\bigg )M \\
= &\ \bigg (\sum_{i=0}^{n} \gamma_ih_i \bigg) \bigg
[\sum_{i=0}^{r} (-1)^i \binom{r}{i} \bigg ( \big (x \otimes
t^{\underline{j} + \underline{k} + (1-i) \underline{h}} \big) \big
(t^{\underline{p} - (1-i) \underline{h}}  D(\underline{\gamma}) \big ) \\
&\ + \big (2x \otimes t^{\underline{j} + \underline{k} - i \underline{h}}
\big ) \big (t^{\underline{p} + i \underline{h}} D(\underline{\gamma})
\big ) +  \big (x \otimes t^{\underline{j} + \underline{k} - (1+i)
\underline{h}} \big ) \big (t^{\underline{p} + (1+i) \underline{h}}
D(\underline{\gamma}) \big ) \bigg) \bigg ]M \\
= &\ \bigg (\sum_{i=0}^{n} \gamma_ih_i \bigg) \bigg
[\sum_{i=0}^{r+2} (-1)^i \binom{r+2}{i} \bigg ( \big (x \otimes
t^{(\underline{j} + \underline{k} + \underline{h}) - i\underline{h}}
\big) \big (t^{(\underline{p} -  \underline{h}) + i \underline{h}}
D(\underline{\gamma}) \big ) \bigg) \bigg ]M
\end{align*}
for all $x \in \mathfrak{g}, \ \underline{j}, \underline{k} \in
\mathbb{Z}^{n+1}, \ \underline{h} \in \mathbb{Z}^{n+1} \setminus
\{\underline{0}\}$ and $\underline{p} \in
\mathbb{Z}\underline{h}$. But since $\underline{\gamma}$ is a generic
vector in $\mathbb{C}^{n+1}$, we can simply take $m =r+2$ and obtain the
desired conclusion.

\noindent (2) If $L(\mathfrak{g})$ acts trivially on $M$, then the
assertion follows immediately from \cite[Theorem 4.8]{BF}. So we assume
that  $L(\mathfrak{g})$ acts non-trivially on $M$. Now since
$\widehat{M}$ is an $A$-module (by Proposition \ref{Cover}), it suffices
to show that one of the weight spaces of $\widehat{M}$ is
finite-dimensional. To this end, first fix any arbitrary $\lambda +
\underline{s} \in P(M)$ for some $\lambda \in D^*$ and $\underline{s} \in
\mathbb{Z}^{n+1}$. Then we have
\begin{align*}
d_i \big (\mu_{x \otimes t^{\underline{s} - \underline{k}},u} \big) = (\lambda + s_i) \big (\mu_{x \otimes t^{\underline{s} - \underline{k}},u} \big) 
\end{align*}
for all $u \in M_{\lambda + \underline{k}}, \ x \in \mathfrak{g}, \ \underline{k} \in \mathbb{Z}^{n+1},\ 0 \leqslant i \leqslant n$, which thereby implies that
\begin{align*}
\widehat{M}_{\lambda + \underline{s}} = \text{span} \{\mu_{x \otimes t^{\underline{s} - \underline{k}},u} \ | \ x \in \mathfrak{g}, \ u \in M_{\lambda + \underline{k}}, \ \underline{k} \in \mathbb{Z}^{n+1} \}. 
\end{align*}
Consider the usual $||\cdot||_1$ norm on $\mathbb{C}^{n+1}$ which is defined by setting\\ $||\underline{w}||_1 = \sum_{i=0}^{n} |w_i| \ \forall \ \underline{w} \in \mathbb{C}^{n+1}$. Now by (1), there exists $m \in \mathbb{N}$ such that 
\begin{align}\label{Annihilators}
\big (\mathcal{T}_{\underline{j}, \underline{0}}^{(m, e_i)}(x) \big ) M = (0) \ \forall \ x \in \mathfrak{g}, \ \underline{j} \in \mathbb{Z}^{n+1}, \ 0 \leqslant i \leqslant n.
\end{align}
\textbf{Claim.} $\widehat{M}_{\lambda + \underline{s}} = \text{span} \bigg \{\mu_{x \otimes t^{\underline{s} - \underline{k}},u} \ \big | \ x \in \mathfrak{g}, \  u \in M_{\lambda + \underline{k}}, \ ||k||_1 \leqslant \dfrac{(n+1)m}{2} \bigg \}$.\\
From Remark \ref{Weights}, we have $M_{\lambda + \underline{k}} = M_ {\underline{\alpha}}$, where $\underline{\alpha} = (\underline{\gamma}, \underline{\beta} + \underline{k})$ and $\underline{\beta} = (\lambda(d_0), \ldots, \lambda(d_n)) \in \mathbb{C}^{n+1}$. Without loss of generality, assume that
\begin{align} \label{Non-zero}
\underline{\beta} = \underline{0} \, \ \text{if} \, \ \underline{\beta} + \mathbb{Z}^{n+1} = \mathbb{Z}^{n+1}, 	
\end{align}
which ensures us that $D(\underline{\gamma})$ acts by a \textit{non-zero}
scalar on every weight space of $M$ except $M_{\underline{\beta}}$. We
proceed to prove our claim by induction on $||k||_1$.\\
If $|k_j| \leqslant m/2 \ \forall \ 0 \leqslant j \leqslant n$, then
$||k||_1 \leqslant \dfrac{(n+1)m}{2}$ and the claim is immediate.
So suppose $|k_j| > m/2$ for some $0 \leqslant j \leqslant n$; in
particular, $\underline{k} \neq \underline{0}$. It suffices to
consider $k_j < -m/2$ (the case $k_j > -m/2$ is similar). Under
this assumption, $||k + ie_j||_1 < ||k||_1 \ \forall \
1 \leqslant i \leqslant m$. Now by (\ref{Non-zero}), for any $u
\in M_{\lambda + \underline{k}}$, there exists $v \in M_{\lambda +
\underline{k}}$ such that $D(\underline{\gamma})v = u$. Applying
(\ref{Annihilators}), we obtain
\begin{align*}
\sum_{i=0}^{m} (-1)^i \binom{m}{i} \mu_{x \otimes t^{\underline{s} -
\underline{k} -ie_j},
(t^{ie_j}D(\underline{\gamma}))v}(t^{\underline{r}})
= &\ \sum_{i=0}^{m} (-1)^i \binom{m}{i} \big (t^{\underline{r}} (x
\otimes t^{\underline{s} - \underline{k} -ie_j}) \big )
(t^{ie_j}D(\underline{\gamma}))v \\
= &\ (\mathcal{T}_{\underline{r} + \underline{s} - \underline{k},
\underline{0}}^{(m, e_j)}(x) \big )v = 0 \ \forall \ \underline{r} \in
\mathbb{Z}^{n+1},
\end{align*}
which thereby yields that
\begin{align}\label{Induction}
\mu_{x \otimes t^{\underline{s} - \underline{k}}, D(\underline{\gamma})v} = - \sum_{i=1}^{m} 
(-1)^i \binom{m}{i} \mu_{x \otimes t^{\underline{s} - (\underline{k} + ie_j)}, (t^{ie_j}D(\underline{\gamma}))v}
\end{align} 
But as $(t^{ie_j}D(\underline{\gamma}))v \in M_{\lambda + (\underline{k} + ie_j)}, \ ||\underline{k} + ie_j||_1 < ||k||_1  \ \forall \ 1 \leqslant i \leqslant m$, we can use induction to see that 
\begin{align*}
\text{the RHS of (\ref{Induction}) belongs to} \  \text{span} \bigg \{\mu_{x \otimes t^{\underline{s} - \underline{k}}, u} \ \big | \ x\in \mathfrak{g}, \ u \in M_{\lambda + \underline{k}}, \ ||k||_1 \leqslant \dfrac{(n+1)m}{2} \bigg \}.
\end{align*}
Hence the claim follows and the lemma is proved.                
\end{proof}

\blmma \label{Quotient}
Let $M$ be an irreducible cuspidal module over $\widehat{\tau}$. Then $M$ is isomorphic to an irreducible quotient of $V(c, \lambda_1, \lambda_2, \underline{\alpha})$ for some quadruplet $(c, \lambda_1, \lambda_2, \underline{\alpha}) \in \mathbb{C} \times P_{\mathfrak{g}}^{+} \times P_{\mathfrak{sl}_{n+1}}^{+} \times \mathbb{C}^{n+1}$
(see Remark \ref{Quadruplet}).  
\elmma
\begin{proof} 
If $L(\mathfrak{g})$ acts trivially on $M$, then $M$ is an irreducible cuspidal module over $\mathcal{W}_{n+1}$ and thus the theorem follows directly from \cite[Lemma 5.5]{BF}. So suppose that $L(\mathfrak{g})$ acts non-trivially on $M$.\\ 
Now since $L(\mathfrak{g})$ is an ideal of $\widehat{\tau}$,  $L(\mathfrak{g})M$ is a non-zero $\widehat{\tau}$-submodule of $M$, which in turn gives $M = L(\mathfrak{g})M$, due to the irreducibility of $M$. Then by Proposition \ref{Cover}, there exists a $\widehat{\tau}$-module homomorphism $\pi : \widehat{M} \longrightarrow M$ satisfying $\pi(\widehat{M}) = M$. Again by Lemma \ref{Cuspidal Cover}, $\widehat{M}$ is a cuspidal module over $A\widehat{\tau}$, which implies that all the weight spaces of $\widehat{M}$ must have the same dimension. Consequently $\widehat{M}$ has finite length, say 
\begin{align*}
(0) = \widehat{M_0} \subseteq \widehat{M_1} \subseteq \ldots \subseteq \widehat{M_{l-1}} \subseteq \widehat{M_l} = \widehat{M},
\end{align*} 
for some $l \in \mathbb{N}$, where $\widehat{M_i}/\widehat{M_{i-1}}$ is irreducible over $A\widehat{\tau}$ for every $1 \leqslant i \leqslant l$.\\
Let $p$ be the smallest integer such that $\pi(\widehat{M_p}) \neq (0)$, which implies that $\pi(\widehat{M_p}) = M$ together with $\pi(\widehat{M_{p-1}}) = (0)$. This induces an onto $\widehat{\tau}$-module homomorphism $\overline{\pi} : \widehat{M_p}/\widehat{M_{p-1}} \longrightarrow M$. Subsequently by Theorem \ref{Jets}, $\widehat{M_p}/\widehat{M_{p-1}} \cong V(c, \lambda_1, \lambda_2, \underline{\alpha})$ for some $(c, \lambda_1, \lambda_2, \underline{\alpha}) \in \mathbb{C} \times (P_{\mathfrak{g}}^{+})^{\times} \times P_{\mathfrak{sl}_{n+1}}^{+} \times \mathbb{C}^{n+1}$, where $(P_{\mathfrak{g}}^{+})^{\times} = (P_{\mathfrak{g}}^{+}) \setminus \{0\}$. But again by Proposition \ref{Fundamental},  $V(c, \lambda_1, \lambda_2, \underline{\alpha})$ is irreducible over $\widehat{\tau}$ as $\lambda_1 \neq 0$ and thus $M \cong V(c, \lambda_1, \lambda_2, \underline{\alpha})$ as $\widehat{\tau}$-modules in this case.   
\end{proof}

We are now ready to state the main theorem of this section which is an immediate consequence of Theorem \ref{T2.13}, Proposition \ref{Fundamental} and Lemma \ref{Quotient} (also see \cite[Theorem 5.4]{BF}). 

\bthm \label{T2.14}
Any non-trivial irreducible cuspidal $\tau$-module is isomorphic
to one of the following.
\begin{enumerate}
	\item $V(c, \lambda_1, \lambda_2, \underline{\alpha})$ for some $(c, \lambda_1, \lambda_2, \underline{\alpha}) \in \mathbb{C} \times (P_{\mathfrak{g}}^{+})^{\times} \times P_{\mathfrak{sl}_{n+1}}^{+} \times \mathbb{C}^{n+1}$, where $(P_{\mathfrak{g}}^{+})^{\times} = (P_{\mathfrak{g}}^{+}) \setminus \{0\}$.
	\item A $\mathcal{W}_{n+1}$-module $V(c, 0, \lambda_2, \underline{\alpha})$ for some $(c, \lambda_2, \underline{\alpha}) \in \mathbb{C} \times P_{\mathfrak{sl}_{n+1}}^{+} \times \mathbb{C}^{n+1}$, where $(\lambda_2, c) \neq (\omega_k,k)$ for any fundamental weight $\omega_k$ of $\mathfrak{sl}_{n+1}$, with $1 \leqslant k \leqslant n$. 
	\item A $\mathcal{W}_{n+1}$-module $V(n+1, 0, 0, \underline{\alpha})$, where $\underline{\alpha} \notin \mathbb{Z}^{n+1}$. 
	\item A $\mathcal{W}_{n+1}$-submodule $\pi_k \big(V(k , 0, \omega_k, \underline{\alpha}) \big)$ of $V(k+1 , 0, \omega_{k+1}, \underline{\alpha})$, where $\underline{\alpha} \in \mathbb{C}^{n+1}$, $0 \leqslant k \leqslant n$ and $\omega_0 = \omega_{n+1} = 0$. 
	\end{enumerate}
\ethm

\brmk \label{Dense}
From our Theorem \ref{T2.14}, we can easily deduce that any non-trivial irreducible cuspidal module $V$ over $\tau$ satisfies the following.
\begin{enumerate}
\item Either $P(V) = \lambda + \mathbb{Z}^{n+1}$ for some $\lambda \in D^*$ or $P(V) = \mathbb{Z}^{n+1} \setminus \{\underline{0}\}$.
\item $\dim V_{\mu_1} = \dim V_{\mu_2} \ \forall \ \mu_1, \mu_2 \in P(V) \setminus \{0\}$. 
\end{enumerate}
\ermk

\section{Generalized Highest Weight Modules}\label{GHW modules}

In this section, we introduce the notion of generalized highest weight modules over $\tau$ and show that every irreducible Harish-Chandra module must be either a cuspidal module or a generalized highest weight module. We conclude this section by recording some of the important properties of generalized highest weight modules.

\noindent \textbf{Notations.} If $\underline{k}, \underline{l} \in \mathbb{Z}^{n+1}$ such that $k_i  \geqslant l_i$ for all $0 \leqslant i \leqslant n$, then we say that $\underline{k} \geqslant \underline{l}$. For $p, q \in \mathbb{Z}$, we set $[p,q] = \{m \in \mathbb{Z} \ | \ p \leqslant m \leqslant q \}$ and define $(-\infty, p]$ and $[q, \infty)$ similarly.

\bdfn
$V$ is said to be a generalized highest weight (GHW) module with GHW $\Lambda_0$ if there exist a non-zero $v_{\Lambda_{0}} \in V_{\Lambda_{0}}$ and some $p \in \mathbb{N}$ such that $\tau_{\underline{k}}v_{\Lambda_{0}} = (0) \ \forall \ \underline{k} \geqslant (p,p, \ldots, p,p)$. In this case, $v_{\Lambda_{0}}$ is called a GHW vector. 
\edfn

To proceed further with our classification problem, we need the following lemma related to cuspidal modules, which is absolutely pivotal for our cause.

\blmma \label{L3.2}
Let $V$ be a cuspidal (but not necessarily irreducible) module over $\tau$. Then:
\begin{enumerate}
	\item $V$ has a $non$-$zero$ irreducible $\tau$-submodule.
	\item $V$ has finite length.
	\item $\dim V_{\mu_1} = \dim V_{\mu_2} \ \forall \ \mu_1, \mu_2
	\in P(V) \setminus \{0\}$.
\end{enumerate}
\elmma
\begin{proof}
	(1) If $V$ has a trivial $\tau$-submodule, then the assertion is clear. So assume that $V$ does not have any trivial $\tau$-submodules. Pick any $0 \neq \mu \in P(V)$ and set $M = \bigoplus_{\underline{k} \in \mathbb{Z}^{n+1}} V_{\mu + \underline{k}}$, which is a non-trivial $\tau$-module. Clearly $\mu = \sum_{i=0}^{n}c_i \delta_i$ for some $ \underline{c} = (c_0, \ldots, c_n) \in \mathbb{C}^{n+1}$, where $\{\delta_0, \ldots, \delta_n \}$ is the usual dual basis of $D$.\\
	\textbf{Claim 1.} Every non-zero submodule of $M$ has a common weight, say $\lambda^{\prime}$.\\
	Let $W$ be any non-zero $\tau$-submodule of $M$. By our initial assumption on $V$ and Lemma \ref{L2.2}, $W$ admits a non-trivial irreducible $\tau$-subquotient, say $W^{\prime}/W^{\prime \prime}$. Also note that $P(W^{\prime}/W^{\prime \prime}) \subseteq \mu + \mathbb{Z}^{n+1}$.\\
	\textbf{Case 1.} If $\underline{c} \notin \mathbb{Z}^{n+1}$, then by Remark \ref{Dense}, $P(W^{\prime}/W^{\prime \prime}) = \mu + \mathbb{Z}^{n+1}$, which thereby implies that $P(W) = \mu + \mathbb{Z}^{n+1}$ and thus the claim follows.\\
	\textbf{Case 2.} For $\underline{c} \in \mathbb{Z}^{n+1}$, we can again use Remark \ref{Dense} and obtain the following subcases.\\
	(i) $P(W^{\prime}/W^{\prime \prime}) = \mathbb{Z}^{n+1}$, whence it follows that $P(W) = \mathbb{Z}^{n+1}$.\\
	(ii) $P(W^{\prime}/W^{\prime \prime}) = \mathbb{Z}^{n+1} \setminus \{\underline{0}\}$, which implies that $P(W) \supseteq \mathbb{Z}^{n+1} \setminus \{\underline{0}\}$.\\
	From (i) and (ii), it is thus evident that every non-zero element of $\mathbb{Z}^{n+1}$ is a weight of $W$ in this case. This proves the first claim. \\
	\textbf{Claim 2.} $U(\tau)N_{\lambda^{\prime}}$ is a non-zero irreducible $\tau$-submodule of $V$, where $N$ is a non-zero $\tau$-submodule of $M$ such that $\dim
	N_{\lambda^{\prime}}$ is minimal. \\ 
	Pick any non-zero vector $v \in U(\tau)N_{\lambda^{\prime}}$ and consider
	the weight space $(U(\tau)v)_{\lambda^{\prime}}$ of the
	$\tau$-submodule $U(\tau)v$ of $N$. Then using the minimality of
	$\dim N_{\lambda^{\prime}}$, we have $(U(\tau)v)_{\lambda^{\prime}} = N_{\lambda^{\prime}}$. This implies that $U(\tau)N_{\lambda^{\prime}} \subseteq U(\tau)v$ and hence the second claim is established.\\
	(2) First apply (1) to obtain a non-zero irreducible submodule of $V$, say $V_1$. Next consider the cuspidal module $V/V_1$ and again apply (1) to obtain a non-zero irreducible submodule of $V/V_1$, say $V_2/V_1$. This induces a chain of cuspidal submodules of $V$, namely $(0) \subsetneq V_1 \subsetneq V_2 \subsetneq V$. But this process must  terminate after finitely many steps, due to Remark \ref{Dense}. Therefore we ultimately get a finite chain of submodules of $V$, given by
	\begin{align*}
		(0) = V_0 \subsetneq V_1 \subsetneq V_2 \subsetneq \ldots \ldots \subsetneq V_{m-1} \subsetneq V_m = V,
	\end{align*}
	where $V_{i+1}/V_i$ is irreducible over $\tau$ for each $0 \leqslant i \leqslant m-1$, with $m \in \mathbb{N}$.\\
	(3) Follows by just combining Remark \ref{Dense} and (2).  
\end{proof}

\bthm\label{T2.9}
Let $V$ be a non-trivial irreducible Harish-Chandra module over $\tau$. Then $V$ is either a cuspidal module or a GHW module up to a change of co-ordinates.
\ethm
\begin{proof} Suppose that $V$ is not a cuspidal module over $\tau$. Recall that $P(V) \subseteq \{\lambda + \underline{m} \ | \ \underline{m} \in \mathbb{Z}^{n+1} \}$ for some $\lambda \in D^*$. Now for each $j \in \mathbb{Z}$, set $M(j) = \bigoplus_{\underline{k} \in \mathbb{Z}^{n+1}, \ k_0 = 0} V_{\lambda +je_0 + \underline{k}}$. Clearly every $M(j)$ is a Harish-Chandra module with respect to $D^{\prime}=$span$\{d_i \ | \ 1 \leqslant i \leqslant n \}$ over the full toroidal Lie algebra $\tau^{\prime}$ in the $n$ variables $t_1, \ldots, t_n$. Again take $M = \bigoplus_{\underline{k} \in \mathbb{Z}^{n+1}, \ k_1 = 0} V_{\lambda + \underline{k}}$, which is also a Harish-Chandra module over the full toroidal Lie algebra $\tau^{\prime \prime}$ in the $n$ variables $t_0, t_2, \ldots, t_n$ with respect to $D^{\prime \prime}= \text{span} \{d_i \ | \ 0 \leqslant i \leqslant n, i \neq 1 \}$.\\
	\textbf{Claim.} Either $M$ is not cuspidal or $M(j)$ is not cuspidal for some $j \in \mathbb{Z}$.\\
	If not, then both $M$ as well as $M(j)$ is cuspidal for each $j
	\in \mathbb{Z}$. So we can find $N_1 \in \mathbb{N}$ such that
	$\dim V_{\lambda + \underline{k}} \leqslant N_1 \ \forall \
	\underline{k} \in \mathbb{Z}^{n+1}$, with $k_1=0$. Putting $N_2 =
	\max \{N_1, \dim V_0 \}$, an immediate application of
	Lemma \ref{L3.2} yields $\dim V_{\lambda + je_0 + \underline{k}}
	\leqslant \max \{\dim V_{\lambda + je_0 + \underline{k} -
	k_1e_1}, \dim V_0 \} \leqslant N_2 \ \forall \ j \in \mathbb{Z}$
	and $\underline{k} \in \mathbb{Z}^{n+1}$, with $k_0=0$. This
	implies that $V$ is a cuspidal module, which is a contradiction. Hence the claim.
	
	Without loss of generality, assume that $M(0)$ is not cuspidal. Therefore there exists some $\underline{k} \in \mathbb{Z}^{n+1}$, with $k_0 = 0$, such that 
	\begin{align}\label{(2.1)}
		\dim V_{\lambda - \underline{k}} > (\dim \mathfrak{g} +
		2n+2) \bigg (\dim V_{\lambda + e_0} + \sum_{i=1}^{n}\dim
		V_{\lambda + e_0 + e_i} \bigg).
	\end{align}
	Set $\underline{s_0} = \underline{k} + e_0, \ 
	\underline{s_i} =  \underline{k} + e_0 + e_i \ \forall \ 1 \leqslant i \leqslant n$. Then the linear transformation sending each $\underline{e_i}$ to $\underline{s_i}$, for $0 \leqslant i \leqslant n$, has determinant $1$ and so $\{ \underline{s_0}, \ldots, \underline{s_n} \}$ is also a $\mathbb{Z}$-basis of $\mathbb{Z}^{n+1}$. Moreover note that $\underline{s_0} + (\lambda - \underline{k}) = \lambda + e_0$ and $\underline{s_i} + (\lambda - \underline{k}) = \lambda + e_0 + e_i \ \forall \ 1 \leqslant i \leqslant n$. Consequently using (\ref{(2.1)}), we can conclude that, up to a change of co-ordinates, there exists a non-zero vector $v_0 \in V_{\lambda - \underline{k}}$ such that $\tau_{\underline{e_i}}v_{0} = (0) \ \forall \ 0 \leqslant i \leqslant n$. Henceforth repeated applications of the commutator relations yield that $V$ is a GHW module with a GHW vector $v_0$ of weight $\lambda - \underline{k}$.

\end{proof}

\blmma\label{L3.5}
Let $V$ be a non-trivial irreducible GHW module over $\tau$. Then:
\begin{enumerate}
	\item Every non-zero $v \in V$ is a GHW vector of $V$. 
	\item $\tau_{-\underline{k}}v \neq (0)$ for any $0 \neq v \in V$ and $\underline{k} \in \mathbb{N}^{n+1}$. 
	\item For each $\mu \in P(V)$ and $\underline{k} \in \mathbb{N}^{n+1}$, $\{m \in \mathbb{Z} \ | \ \mu + m \underline{k} \in P(V) \} = (-\infty, l]$ for some $l \in \mathbb{Z}_{+}$. 
\end{enumerate}
\elmma
\begin{proof}
	The proofs of (1) and (2) are similar to \cite[Lemma 3.3]{MZ} and \cite[Lemma 3.4]{MZ} respectively.\\
	(3) 	Set $J = \{m \in \mathbb{Z} \ | \ \mu + m \underline{k} \in P(V) \}$. 
	From Lemma \ref{L3.5}, we have either $J = (-\infty, l]$ for some $l \in \mathbb{Z_{+}}$ or $J = \mathbb{Z}$. We claim that $J \neq \mathbb{Z}$. If not, then we must have $J = \mathbb{Z}$. To this end, define
	\begin{align*}
	L_{\underline{k}} := \text{span} \{t^{r\underline{k}}d_i, t^{s\underline{k}}K_j \ | \ r, s \in \mathbb{Z}, \ 0 \leqslant i, j \leqslant n \}.
	\end{align*}
	Putting $M_{\underline{k}} = \bigoplus_{p \in \mathbb{Z}} V_{\mu +
		p\underline{k}}$, 
	 the proof now proceeds verbatim as in \cite[Lemma 3.5]{MZ}.
\end{proof}

\blmma \label{L3.8}
Let $V$ be an irreducible Harish-Chandra module over $\tau$ which is not cuspidal. Then up to a possible change of co-ordinates, $V$ satisfies the following properties.
\begin{enumerate}
	\item $V$ is a GHW module over $\tau$ having a GHW vector $v_0$ of weight $\Lambda_0$. 
	\item $\Lambda_{0} + \underline{k} \notin P(V) \ \forall \ \underline{0} \neq \underline{k} \in \mathbb{Z}_{+}^{n+1}$. 
	\item $\Lambda_{0} - \underline{k} \in P(V) \ \forall \   \underline{k} \in \mathbb{Z}_{+}^{n+1}$. 
	\item For any $\underline{k}, \underline{l} \in \mathbb{Z}^{n+1}$ with $\underline{k} \leqslant \underline{l}$, the condition $\Lambda_{0} + \underline{k} \notin P(V)$ implies $\Lambda_{0} + \underline{l} \notin P(V)$. 
	\item For any $\underline{0} \neq \underline{k} \in \mathbb{Z}_{+}^{n+1}$ and $\underline{l} \in \mathbb{Z}^{n+1}$, we have\\
	$\{m \in \mathbb{Z} \ | \ \Lambda_{0} +  \underline{l} + m \underline{k} \in P(V) \} = (-\infty, q]$ for some $q \in \mathbb{Z}$.
\end{enumerate} 
\elmma  
\begin{proof} By Theorem \ref{T2.9}, $V$ is a non-trivial GHW module. Consequently (1), (2), (3) and (4) follow verbatim as in \cite[Lemma 3.6]{MZ}.\\
	(5) This can be deduced by considering the $\mathbb{Z}$-basis $\{e_i^{\prime}\}_{i=0}^{n}$ of $\mathbb{Z}^{n+1}$ as in \cite[Lemma 3.6]{MZ} along with Lemma \ref{L3.5} and finally performing a suitable change of co-ordinates.    	
\end{proof} 

\section{Modules of the Highest Weight Type}\label{Highest Weight}
In this section, we define irreducible highest weight modules for $\tau$ with respect to a suitable triangular decomposition and then provide concrete realizations of their highest weight spaces.

Consider a $\mathbb{Z}$-grading on $\tau$ by means of eigenvalues of
$d_0$. This naturally gives rise to a triangular decomposition of $\tau$,
namely
\begin{align*}
	\tau = &\ \tau_{-} \ \oplus \ \tau_{0} \ \oplus \ \tau_{+}, \
	\text{where}\\
	\tau_{+} = &\ \mathfrak{g} \otimes t_0
	\mathbb{C}[t_0, t_1^{\pm 1}, \ldots, t_n^{\pm 1}]  
	 \ \oplus \ \sum_{i=0}^{n} t_0 \mathbb{C}[t_0, t_1^{\pm 1},
	\ldots, t_n^{\pm 1}]K_i \ \oplus \ \sum_{i=0}^{n} t_0
	\mathbb{C}[t_0, t_1^{\pm 1}, \ldots, t_n^{\pm 1}]d_i,\\
	\tau_{0} = &\ \mathfrak{g} \otimes \mathbb{C}[t_1^{\pm 1},
	\ldots, t_n^{\pm n}] \ \oplus \ \sum_{i=0}^{n} \mathbb{C}[t_1^{\pm 1}, \ldots, t_n^{\pm 1}]K_i \ \oplus \ \sum_{i=0}^{n} \mathbb{C}[t_1^{\pm 1}, \ldots, t_n^{\pm 1}]d_i,\\
	\tau_{-} = &\ \mathfrak{g} \otimes t_0^{-1}
	\mathbb{C}[t_0^{-1}, t_1^{\pm 1}, \ldots, t_n^{\pm 1}] 
	 \ \oplus \ \sum_{i=0}^{n} t_0^{-1} \mathbb{C}[t_0^{-1},
	t_1^{\pm 1}, \ldots, t_n^{\pm 1}]K_i \ \oplus \ \sum_{i=0}^{n}
	t_0^{-1} \mathbb{C}[t_0^{-1}, t_1^{\pm 1}, \ldots, t_n^{\pm
	1}]d_i. 
\end{align*}

The above decomposition enables us to define irreducible modules of the highest weight type for $\tau$. Let $X$ be an irreducible   
$\tau_{0}$-module. Postulate $\tau_{+}X = (0)$, which in turn gives rise to a $(\tau_{0} \oplus \tau_{+})$-module structure on $X$. Define the generalized Verma module $M(X)$ over $\tau$  by setting
\begin{align*}
	M(X) = U(\tau) \otimes_{U(\tau_{0} \oplus \tau_{+})} X.
\end{align*}     
By standard arguments, it can be shown that $M(X)$ has a unique irreducible quotient, which we shall denote by $L(X)$. Again as $\tau$ is a $\mathbb{Z}^n$-extragraded Lie algebra, it follows from \cite[Theorem 1.5]{BZ} (also see \cite[Theorem 1.12]{BB}) that $L(X)$ is a \textit{Harish-Chandra} module, provided that $X$ is a \textit{cuspidal} $\mathbb{Z}^n$-graded \textit{exp-polynomial} module (refer to \cite[Section 1]{BZ} for more details on $\mathbb{Z}^n$-extragraded Lie algebras and $\mathbb{Z}^n$-graded exp-polynomial modules). Finally $L(X)$ is a cuspidal $\tau$-module if and only if $X$ is the trivial module (in which case $L(X)$ itself is the trivial module over $\tau$).

Now there exists an action of $GL(n+1, \mathbb{Z})$ (by automorphisms) on $\tau$ (see Subsection \ref{Co-ordinates}). So for each $A \in GL(n+1, \mathbb{Z})$, we obtain a twisted irreducible module over $\tau$, which we shall denote by $L(X)^A$, where the action of $\tau$ is twisted by an automorphism $A$ of $\mathbb{Z}^{n+1}$. We shall refer to these twisted modules of the form $L(X)^A$ as simply \textit{irreducible modules of the highest weight type}. 	
\brmk \label{Yuly} 
The classification of irreducible $\tau$-modules of the highest weight type of \textit{non-zero level} (i.e. one of the $K_i$'s acts \textit{non-trivially} on the module) has been carried out in \cite{B}. In this instance, the author refers to these modules as \textit{bounded modules}, in the sense that the real parts of eigenvalues (with respect to $d_0$) on these modules are bounded from above. Towards the conclusion of this paper, we shall also classify all \textit{level zero} (i.e. each $K_i$ acts trivially on the module) irreducible bounded modules, thereby finally obtaining all the irreducible bounded modules up to isomorphism.
\ermk

\bthm \label{L(X)}
Let $L(X)$ be an irreducible Harish-Chandra module over $\tau$ of the highest weight type, where $X$ is an irreducible cuspidal $\tau_0$-module. 
\begin{enumerate}
\item If both $K_0$ and $d_0$ act trivially on $X$, then $t^{\underline{m}}K_0$ as well as  $t^{\underline{m}}d_0$ act trivially on $X$ for all $\underline{m} \in \{\underline{0}\} \times \mathbb{Z}^n$ and consequently $X$ is an irreducible cuspidal module over the full toroidal Lie algebra in the $n$ variables $t_1, \ldots, t_n$, that we have already classified in Theorem \ref{T2.14}.
	
\item If either $K_0$ or $d_0$ acts non-trivially on $X$, then there exist a finite-dimensional irreducible $\mathfrak{g}$-module $V_1$, a finite-dimensional irreducible $\mathfrak{gl}_n$-module $V_2$,  $\underline{\alpha} = (\alpha_1, \ldots, \alpha_n) \in \mathbb{C}^n$ and  $(a, b) \in \mathbb{C}^2 \setminus \{(0,0)\}$ such that   
\begin{align*}
	X \cong V_1 \otimes V_2 \otimes \mathbb{C}[t_1^{\pm 1}, \ldots, t_n^{\pm n}],
\end{align*}	
where the action of $\tau_0$ on $X$ is given by 
\begin{align*}
	(x \otimes t^{\underline{m}}) (v_1 \otimes v_2 \otimes
	t^{\underline{r}}) = &\ (xv_1) \otimes v_2 \otimes t^{\underline{r}
	+ \underline{m}}, \\ 
	t^{\underline{m}}K_0 (v_1 \otimes v_2 \otimes t^{\underline{r}})
	= &\ a(v_1 \otimes v_2 \otimes t^{\underline{r} +
	\underline{m}}), \\
	t^{\underline{m}}K_i (v_1 \otimes v_2 \otimes t^{\underline{r}})
	= &\ 0, \\
	t^{\underline{m}}d_0 (v_1 \otimes v_2 \otimes t^{\underline{r}})
	= &\ b(v_1 \otimes v_2 \otimes t^{\underline{r} +
	\underline{m}}), \\
	t^{\underline{m}}d_i (v_1 \otimes v_2 \otimes t^{\underline{r}})
	= &\ (\alpha_i + r_i)(v_1 \otimes v_2 \otimes t^{\underline{r} +
	\underline{m}}) + \sum_{j=1}^{n}m_j \big (v_1 \otimes
	(E_{j,i}v_2) \otimes t^{\underline{r} + \underline{m}} \big ) \\
	&\ \forall \ x \in \mathfrak{g}, \ v_1 \in V_1, \ v_2 \in V_2, \
	\underline{m},\underline{r} \in \{\underline{0}\} \times \mathbb{Z}^n, \ 1 \leqslant i
	\leqslant n.    
\end{align*}
Here $E_{j,i}$ denotes the matrix of order $n$ having $1$ at the $(j,i)$-th entry and $0$ elsewhere.
\end{enumerate}
\ethm
\begin{proof}
Note that both $K_0$ and $d_0$ are central elements of $\tau_0$ and so they act by fixed scalars, say $a$ and $b$ respectively on $X$. 
Now $X$ is a cuspidal module over the full toroidal Lie algebra $\tau^{\prime}$ in the $n$ variables $t_1, \ldots, t_n$, whence by Lemma \ref{L3.2}, we obtain an irreducible $\tau^{\prime}$-submodule of $X$, say $X_0$. Then by Theorem \ref{T2.13}, $\mathcal{Z}^{\prime} = \text{span} \{t^{\underline{m}}K_i \ | \ \underline{m} \in \{\underline{0}\} \times \mathbb{Z}^n,\ 1 \leqslant i \leqslant n \}$ acts trivially on $X_0$. But $W = \text{span} \{w \in X \ | \ \mathcal{Z}^{\prime}w = (0) \}$ is a $\tau_0$-submodule of $X$, which implies that $\mathcal{Z}^{\prime}$ acts trivially on $X$, thanks to the irreducibility of $X$. In particular, $K_1, \ldots, K_n$ act trivially on $X$.\\
(1) $X$ is cuspidal over $(\mathcal{V}ir)_{0}^{\prime \prime} = \text{span} \{t^{\underline{r}}K_0, t^{\underline{s}}d_j \ | \ \underline{r}, \underline{s} \in \{\underline{0}\} \times \mathbb{Z}^n, \ 1 \leqslant j \leqslant n \}$ (and also over $(\mathcal{V}ir)_{0}^{\prime} = \text{span} \{t^{\underline{r}}d_j  \ | \ \underline{r} \in \{\underline{0}\} \times \mathbb{Z}^n, \ 0 \leqslant j \leqslant n \}$), with the  Witt algebra $\mathcal{W}_n$ embedded as a subalgebra inside $(\mathcal{V}ir)_{0}^{\prime \prime}$ (respectively $(\mathcal{V}ir)_{0}^{\prime}$). Hence \cite[Lemma 3.3]{Ma} and Remark \ref{R1.1} reveals that $X$ has finite length over $(\mathcal{V}ir)_{0}^{\prime \prime}$ (and over $(\mathcal{V}ir)_{0}^{\prime}$), which induces an irreducible cuspidal $(\mathcal{V}ir)_{0}^{\prime \prime}$-submodule (respectively $(\mathcal{V}ir)_{0}^{\prime}$-submodule) of $X$, which we denote by $X^{\prime \prime}$ (respectively by $X^{\prime}$).\\ 
\textbf{Case 1.} Let us first take $n=1$.\\ 
In this case, both $(\mathcal{V}ir)_{0}^{\prime}$  as well as $(\mathcal{V}ir)_{0}^{\prime \prime}$ can be identified with $\text{HVir}/ \text{span} \{C_D, C_{DI}, C_I \}$ (see Subsection \ref{Twisted}). Now since $K_0=d_0=0$ on $X$, a straightforward application of \cite[Theorem 4.4]{LZ2} shows that $t_1^{m}K_0$ acts trivially on $X^{\prime \prime}$ and $t_1^{m}d_0$ acts trivially on $X^{\prime}$ for all $m \in \mathbb{Z}$. Finally note that
\begin{align*}
W^{\prime} = \{w \in X \ | \ (t_1^{m}d_0)w = 0 \ \forall \ m \in \mathbb{Z} \} \ \text{and} \
 W^{\prime \prime} = \{w \in X \ | \ (t_1^{m}K_0)w = 0 \ \forall \ m \in \mathbb{Z} \} 
\end{align*}  
are both $\tau_0$-submodules of $X$. Thus we are done with this case, thanks to the irreducibility of $X$.\\
\textbf{Case 2.} Next let us proceed to the case $n \geqslant 2$.\\
In this instance, we can directly apply \cite[Theorem 3.3]{GLZ} to both $(\mathcal{V}ir)_{0}^{\prime}$ and $(\mathcal{V}ir)_{0}^{\prime \prime}$ to infer that $t^{\underline{m}}K_0$ acts trivially on $X^{\prime \prime}$ and $t^{\underline{m}}d_0$ acts trivially on $X^{\prime}$ for all $\underline{m} \in \{\underline{0}\} \times \mathbb{Z}^n$, since $K_0$ and $d_0$ act trivially on $X$. Then analogous to the first case, we can easily deduce that both $t^{\underline{m}}K_0$ and $t^{\underline{m}}d_0$ act trivially on $X$, using the irreducibility of $X$ over $\tau_0$. This concludes the proof of our first assertion. \\
(2) If $K_0$ acts by a non-zero scalar on $X$, then we have $a \neq 0$ and the result follows directly from \cite[Theorem 3.4]{B}. So suppose that $d_0$ acts by a non-zero scalar on $X$. This gives $b \neq 0$, whence a similar argument as in Case 1 of Proposition
\ref{Fundamental} shows that $V_1 \otimes V_2 \otimes
\mathbb{C}[t_1^{\pm 1}, \ldots, t_n^{\pm n}]$ is irreducible  over
$\tau_0^{\prime \prime} = \{x \otimes t^{\underline{r}}, \
t^{\underline{s}}d_i \ | \ x \in \mathfrak{g}, \ \underline{r},
\underline{s} \in \{\underline{0}\} \times \mathbb{Z}^n, \ 0 \leqslant i \leqslant n \}$. Now
$d_0$ plays exactly the same role as $K_0$ on $X$ and
$[t^{\underline{r}}d_0, t^{\underline{s}}K_i] = [t^{\underline{r}}d_0, x
\otimes t^{\underline{s}}] = 0 \ \forall \ x \in \mathfrak{g}, \
\underline{r}, \underline{s} \in \{\underline{0}\} \times \mathbb{Z}^n, \ 1 \leqslant i \leqslant
n$. Hence we can deduce the desired result in an analogous manner
as in \cite[Theorem 3.1]{CD} and \cite[Theorem 3.1]{EJ}, along with an
application of Theorem \ref{Jets}, by just working with $d_0$ instead of
$K_0$.
\end{proof}

\brmk \label{R5.3}
\
\begin{enumerate}
\item From Theorem \ref{L(X)}, it is clear that if $L(X)$ is not the trivial $\tau$-module, then we have either $P(X) = \lambda + \mathbb{Z}^{n}$ for some $\lambda \in D^*$ or $P(X) = \mathbb{Z}^{n} \setminus \{\underline{0}\}$.
\item Note that the $\tau_0$-module $V_1 \otimes V_2 \otimes \mathbb{C}[t_1^{\pm 1}, \ldots, t_n^{\pm 1}]$ in Theorem \ref{L(X)} is fully parametrized by $(a,b,c, \lambda_1, \lambda_2, \underline{\alpha}) \in \mathbb{C} \times \mathbb{C} \times \mathbb{C} \times P_{\mathfrak{g}}^{+} \times P_{\mathfrak{sl}_n}^{+} \times \mathbb{C}^{n}$. So we shall denote this module by $V(a, b, c, \lambda_1, \lambda_2, \underline{\alpha})$, even when $(a,b) =(0,0)$ (equipped with exactly the same $\tau_0$-action as in (2) of Theorem \ref{L(X)} if $(a,b) =(0,0)$). $V(a, b, c, \lambda_1, \lambda_2, \underline{\alpha})$ is irreducible over $\tau_0$ when $(a,b) \neq (0,0)$. But if $(a,b) = (0,0)$, then $V(a, b, c, \lambda_1, \lambda_2, \underline{\alpha}) \cong V(c, \lambda_1, \lambda_2, \underline{\alpha})$ is actually a module over the full toroidal Lie algebra $\tau^{\prime}$ in the $n$ variables $t_1, \ldots, t_n$, which is not necessarily irreducible over $\tau^{\prime}$, but it always admits a unique irreducible quotient. Thus we shall denote the unique irreducible quotient of $V(a, b, c, \lambda_1, \lambda_2, \underline{\alpha})$  by $V^{\prime}(a, b, c, \lambda_1, \lambda_2, \underline{\alpha})$ for every $(a,b,c, \lambda_1, \lambda_2, \underline{\alpha}) \in \mathbb{C} \times \mathbb{C} \times \mathbb{C} \times P_{\mathfrak{g}}^{+} \times P_{\mathfrak{sl}_n}^{+} \times \mathbb{C}^{n}$ (in the framework of the higher-dimensional Virasoro algebra $\mathcal{V}ir$, we shall denote this unique irreducible quotient of $V(a, b, c, \lambda, \underline{\alpha})$  by $V^{\prime}(a, b, c, \lambda, \underline{\alpha})$ for each $(a,b,c, \lambda, \underline{\alpha}) \in \mathbb{C} \times \mathbb{C} \times \mathbb{C} \times  P_{\mathfrak{sl}_n}^{+} \times \mathbb{C}^{n}$).  

\item From the actions of $\tau_0$ and $\tau^{\prime}$ prescribed in Theorem \ref{L(X)}, it can be easily verified that $V^{\prime}(a,b,c, \lambda_1, \lambda_2, \underline{\alpha})$ is a cuspidal $\mathbb{Z}^n$-graded exp-polynomial module over these Lie algebras.
       
\item Towards the end of Section \ref{S6}, we shall show that the cuspidality condition on the irreducible $\tau_0$-module $X$ is redundant if $L(X)$ is a Harish-Chandra module over $\tau$.  
\end{enumerate}
\ermk
\section{Classification of Irreducible GHW Modules} \label{S6}
In this section, unless explicitly stated, $V$ will always denote an irreducible $non$-$trivial$ GHW module over $\tau$, having a GHW vector $v_{\Lambda_{0}}$ of weight $\Lambda_{0}$. This gives $V = \bigoplus_{\underline{k} \in \mathbb{Z}^{n+1}} V_{\Lambda_{0} + \underline{k}}$. Moreover Lemma \ref{L3.5} permits us to take (and we shall take) $\Lambda_0 = 0$ if $0 \in P(V)$. Under these assumptions, our primary aim in this section is to prove the following result. 

\bthm \label{T4.2}
$V \cong L(X)^A$ for some irreducible cuspidal module $X$ over $\tau_0$ and $A \in GL(n+1, \mathbb{Z})$.
\ethm

We require a lot of preparation to prove this theorem. We start  with the following proposition, which will play a crucial role in this regard.

\bppsn \label{P4.3} 
Let $\underline{m}, \ \underline{k} = (k_0, \ldots, k_n) \in \mathbb{Z}^{n+1}$ such that $k_0, \ldots, k_n$ are relatively prime. If there exists $\underline{p} \in \mathbb{N}^{n+1}$ satisfying
\begin{align*}
\bigg \{\Lambda_0 + \sum_{i=0}^{n}m_ie_i + \sum_{i=0}^{n}r_ip_ie_i \ | \
\underline{r} \in \mathbb{Z}^{n+1}, \ \sum_{i=0}^{n}k_ip_ir_i = 0 \bigg
\} \bigcap P(V) = \emptyset,
\end{align*}
then $V \cong L(X)^A$ for an irreducible cuspidal module $X$ over $\tau_0$ and  some $A \in GL(n+1, \mathbb{Z})$.  
\eppsn
\begin{proof} Put $G = \big \{\sum_{i=0}^{n}l_ie_i \in \mathbb{Z}^{n+1} \ | \ \sum_{i=0}^{n}k_il_i = 0 \big \}$. Using Lemma \ref{L3.8} and by performing a change of
	co-ordinates, if necessary, we can show that there exists a
	unique $m_0 \in \mathbb{Z}$ satisfying the following properties
	(see (1) and (2) of \cite[Lemma 3.3]{LZ1}). 
		\begin{align*}
	\big \{\Lambda_{0} + \sum_{i=0}^{n}l_ie_i \in P(V) \ | \ \underline{l} \in \mathbb{Z}^{n+1}, \ \sum_{i=0}^{n}k_il_i \geqslant m_0 \big \} = \phi \ \text{and} \\
	S := \big \{\Lambda_{0} + \sum_{i=0}^{n}l_ie_i \in P(V) \ | \ \underline{l} \in \mathbb{Z}^{n+1}, \ \sum_{i=0}^{n}k_il_i = m_0-1 \big \} \neq \phi.
	\end{align*}	
	Since $k_0, \ldots, k_n$ are relatively prime, there exist $s_0, \ldots, s_n \in \mathbb{Z}$ such that $\sum_{i=0}^{n}k_is_i=1$. Fixing $\underline{t_0} = \sum_{i=0}^{n}s_ie_i$, it is easy to see that
	\begin{align}\label{(5.1)}
		\mathbb{Z}^{n+1}= \mathbb{Z}\underline{t_0} \oplus G \
		\text{and} \ (S + G + \underline{t_0}) \bigcap P(V)
		= \emptyset.
	\end{align}
	Fix $\mu_0 \in S$. It can be readily checked that
	\begin{align}\label{(5.2)}
		S = (\mu_0 + G) \bigcap P(V).
	\end{align}
	Note that the quotient group $\mathbb{Z}^{n+1}/G$ is torsion-free of rank $1$, which permits us to take $G = \mathbb{Z}^n \cong \text{span} \{e_1, \ldots, e_n \}$ and $\underline{t_0} = e_0$, after an appropriate change of co-ordinates. Due to this change of co-ordinates, the weight $\mu_0 \in S$ should also transform accordingly, but for the sake of notational convenience, we shall stick to our notation $\mu_0$ throughout this lemma. Consider $W = \bigoplus_{\lambda \in \mu_0 + \mathbb{Z}^n} V_{\lambda}$, which is a $\tau_0$-submodule of $V$. Pick any $0 \neq v_0 \in V_{\mu_0 + \underline{m^{\prime}}}$, where $\underline{m^{\prime}} \in \{\underline{0}\} \times \mathbb{Z}^n$ is arbitrary.\\
	\textbf{Claim.} $\tau_{-\underline{m^{\prime}}-e_0}v_0 \neq (0)$.\\
	If not, then $\tau_{-\underline{m^{\prime}}-e_0}v_0 =(0)$. Now
	$\mathbb{Z}^{n+1}$ is generated by $\{\underline{m}
	+ e_0,  -\underline{m^{\prime}} - e_0 \ | \ \underline{m} \in
	\{\underline{0}\} \times \mathbb{Z}^n \}$ for any $\underline{m^{\prime}} \in
	\{\underline{0}\} \times \mathbb{Z}^n$. Subsequently $\tau$ is generated, as a Lie
	algebra, by $\tau_{-\underline{m^{\prime}} - e_0}$ and
	$\tau_{\underline{m} + e_0}$, with $\underline{m} \in
	\{\underline{0}\} \times \mathbb{Z}^n$. In view of (\ref{(5.1)}), combined with the
	irreducibility of $V$, this reveals that $\tau$ acts trivially on
	$V$, which contradicts our initial assumption. Hence the claim.
	
	Consequently for every $\underline{m^{\prime}} \in \{\underline{0}\} \times \mathbb{Z}^n$,
	we have an injective linear map
	\begin{align*}
		T_{\underline{m^{\prime}}}: V_{\mu_0 +
		\underline{m^{\prime}}} \longrightarrow (V_{\mu_0 -
		e_0})^{\oplus N}, \ \text{where} \ N= \dim \tau_{-\underline{m^{\prime}} - e_0}.
	\end{align*}     
	Thus $W$ is cuspidal over $\tau_0$ and hence a cuspidal module over the full toroidal Lie algebra $\tau^{\prime}$ in the $n$ variables $t_1,
	\ldots, t_n$ with respect to $D^{\prime} = \text{span} \{d_1,
	\ldots, d_n\}$. By Lemma \ref{L3.2}, $W$ has finite length and
	therefore it admits an irreducible $\tau_0$-submodule, say $X$.
	Then $\tau_{+} X = (0)$ by (\ref{(5.1)}), (\ref{(5.2)}) and Lemma
	\ref{L3.8} and so we are done.    
\end{proof}
We now separately study the case $n=1$ and finally prove Theorem \ref{T4.2} in this particular set-up. So our next few results preceding Theorem \ref{T4.10} will specifically deal with the case $n=1$. 

\blmma \label{L4.4}
Any irreducible Harish-Chandra module over $\tau_{0} = \text{span}\{x \otimes  t_{1}^p, t_{1}^qK_0, t_{1}^rd_0, t_{1}^sd_1, K_1 \ | \ x \in \mathfrak{g}, \ p, q, r, s \in \mathbb{Z} \}$ is either cuspidal or a highest weight module or a lowest weight module. 
\begin{proof} Let $V$ be an irreducible Harish-Chandra module (with respect to $\mathbb{C}d_1$) over $\tau_0$ which is neither a lowest weight module nor cuspidal. Fix $\Lambda \in P(V)$. Then there exists $m \in \mathbb{N}$ such that
	\begin{align*}
		\dim V_{\Lambda-m} > (\dim \mathfrak{g}+3) \
		\dim V_{\Lambda} + 3 \dim V_{\Lambda+1}.  
	\end{align*}
	Thus there exists $0 \neq v \in V_{\Lambda - m}$ such that $(x \otimes t_1^{m})v = t_1^{m}K_0v = t_1^{m+1}K_0v = t_1^{m}d_0v = t_1^{m+1}d_0v = t_1^{m}d_1v = t_1^{m+1}d_1v = 0 \ \forall \ x \in \mathfrak{g}$. Now by repeated use of the commutator relations on $\tau_{0}$, we can deduce that $(\tau_0)_{j}v = (0) \ \forall \ j \geqslant m^2$, where $(\tau_{0})_{j}= \text{span}\{x \otimes  t_{1}^j, t_{1}^jK_0, t_{1}^jd_0, t_{1}^jd_1, \delta_{0,j}K_1 \ | \ x \in \mathfrak{g} \}$. The lemma now directly follows from \cite[Lemma 1.6]{OM}.
\end{proof}
\elmma
\bppsn \label{P4.5}
Let $\Lambda \in P(V)$ such that $(\Lambda + \mathbb{N}e_0 +
\mathbb{Z}e_1) \bigcap P(V) = \emptyset$. Then $V \cong L(X)^A$ for an irreducible cuspidal module $X$ over $\tau_0$ and  some $A \in GL(2, \mathbb{Z})$. 
\eppsn
\begin{proof}
	By hypothesis, we have $V_{+} := \{v \in V \ | \ \tau_{+}v = (0) \} \neq (0)$. Now since $V$ is irreducible, an application of the PBW theorem readily yields that $V_{+}$ is an irreducible Harish-Chandra module with respect to $\mathbb{C}d_1$ over $\tau_0$. Also note that
	\begin{align*}
		V_{+} = \bigoplus_{m \in \mathbb{Z}} V_{+}(\Lambda + m), \ \text{where} \ V_{+}(\Lambda + m) = \{v \in V_{+}\ | \ d_1v = (\Lambda(d_1)+m)v \} \ \forall \ m \in \mathbb{Z}.
	\end{align*}
	\textbf{Claim.} $V_{+}$ is a cuspidal module over $\tau_0$.\\
	If not, then  $V_{+}$ is either a highest weight module or a lowest weight module by Lemma \ref{L4.4}. Suppose that $V_{+}$ is a highest weight module. So there exists $0 \neq v_0 \in V_{+}$ such that $(\tau_{0})_{j}v_{0} = (0) \ \forall \ j \in \mathbb{N}$. By our assumption, $(\tau_0)_{0}$ must act non-trivially on the highest weight space $\{v \in V_{+} \ | \ (\tau_0)^{+}v = (0) \}$ of $V_{+}$, whence there exist a vector $v_{\mu_0}$ of weight $\mu_0$ (say) in $V_{+}$ and $x_0 \in \mathfrak{h} \oplus \sum_{i=0}^{1} \mathbb{C}K_i \oplus \sum_{i=0}^{1} \mathbb{C}d_i$ such that $x_0$ acts by a \textit{non-zero} scalar at $v_{\mu_0}$, where $\mathfrak{h}$ is a Cartan subalgebra of $\mathfrak{g}$. Consequently $V = U(\tau)v_{\mu_0}$ is a highest weight module with respect to the lexicographic order $\prec$ on $\mathbb{Z}^2$ (i.e. $(x_1,x_2) \prec (y_1,y_2) \iff \ \text{either} \ x_1 < y_1 \ \text{or} \ x_1 = y_1 \ \text{and} \ x_2 < y_2$) having highest weight $\mu_0$ and
	\begin{align}\label{P(V)}
	 P(V) \subseteq(\mu_0 - \mathbb{N}e_0 + \mathbb{Z}e_1) \bigcup (\mu_0 - \mathbb{Z_{+}}e_1).
	 \end{align}
	 Thus $V$ is a GHW module having generalized highest weight
	 $\mu_0$, with $(\mu_{0} + \mathbb{N}e_0 + \mathbb{Z}e_1) \bigcap
	 P(V) = \emptyset$ and $(\mu_{0} + \mathbb{Z}e_1) \bigcap
	 P(V) \neq \emptyset$. Hence we can provide a similar
	 argument as in Proposition \ref{P4.3} to conclude that
	 $V^{\prime} = \bigoplus_{r \in \mathbb{Z}} V_{\mu_0 +
	 re_1}$ is a cuspidal module over $\tau_{0}$. Now consider the non-trivial
	 submodule $V^{\prime \prime} = U(\tau_{0})v_{\mu_{0}}$
	 of $V^{\prime}$ and then obtain a
	 non-zero irreducible $\tau_0$-quotient of $V^{\prime \prime}$ (using Zorn's Lemma), say $W^{\prime \prime}$. Again since we already know that
	 $x_0 \in \mathfrak{h} \oplus \sum_{i=0}^{1} \mathbb{C}K_i \oplus
	 \sum_{i=0}^{1} \mathbb{C}d_i$ acts by a non-zero scalar at
	 $v_{\mu_0}$, it is evident that $W^{\prime \prime}$ is a
	 non-trivial cuspidal $\tau_0$-module. Subsequently our Remark
	 \ref{R5.3} yields that $(\mu_0 + \mathbb{Z}e_1) \setminus \{0\}
	 \subseteq P(W^{\prime \prime}) \subseteq P(V)$, which is a
	 contradiction to (\ref{P(V)}). This contradiction establishes
	 our claim and proves the proposition by simply taking $X
	 = V_{+}$.    
\end{proof}
We now record some results in the following lemma, which can be deduced using Lemma \ref{L3.8}, Proposition \ref{P4.3} and Proposition \ref{P4.5}, in a more or less similar manner as in \cite{LZ1,S}.
\blmma \label{L4.6}
In each of the following cases, $V \cong L(X)^A$ for an irreducible cuspidal module $X$ over $\tau_0$ and some $A \in GL(2, \mathbb{Z})$.
\begin{enumerate}
	\item There exist $(i,j) \in \mathbb{Z}^2, \ (k,l) \in \mathbb{Z}^2 \setminus \{(0,0)\}$ and $p, q \in \mathbb{Z}$ such that 
	\begin{align*}
		\{m \in \mathbb{Z} \ | \ \Lambda_{0} + (i,j) + m(k,l) \in P(V) \} \supseteq (-\infty, p] \cup [q, \infty).
	\end{align*}
	\item There exist $(i,j) \in \mathbb{Z}^2$ and $(k,l) \in \mathbb{Z}^2 \setminus \{(0,0)\}$ such that
	\begin{align*}
		\{\Lambda_0 + (i,j) + m(k,l) \ | \ m \in \mathbb{Z} \}
		\bigcap P(V) = \emptyset.
	\end{align*}
		\item There exist $(i,j), (k,l) \in \mathbb{Z}^2$ and $m_1, m_2, m_3 \in \mathbb{Z}$ with $m_1 < m_2 < m_3$ such that 
	\begin{align*}
	& \Lambda_0 + (i,j) + m_1 (k,l) \notin P_D(V), \\
	& \Lambda_0 + (i,j) + m_2(k,l) \in P_D(V), \\
	& \Lambda_0 + (i,j) + m_3(k,l) \notin P_D(V).
	\end{align*}
\end{enumerate}
\elmma


In view of Lemma \ref{L4.6}, it is easy to see that for any $(i,j) \in \mathbb{Z}^2$ and $(k,l) \in \mathbb{Z}^2 \setminus \{(0,0)\}$, there exists some $p \in \mathbb{Z}$ such that
\begin{align}\label{4.3}
\{m \in \mathbb{Z} \ | \ \Lambda_{0} + (i,j) + m(k,l) \in P(V) \} = (-\infty, p] \ \text{or} \ [p, \infty).
\end{align}
Consequently from Lemma \ref{L3.8}, it follows that for each $i \in \mathbb{N}$, there exist $a_i, b_i \in \mathbb{Z}_{+}$ satisfying
\begin{align*}
b_i = &\ \max \{b \in \mathbb{Z} \ | \ \Lambda_0 + (-i,b) \in P(V) \}, \,\,\,\	a_i = \ \text{max} \{a \in \mathbb{Z} \ | \ \Lambda_0 + (a,-i) \in P(V) \}.
\end{align*} 
Then we can deduce the following results from \cite[Claim 1]{LZ1} and \cite[Claim 2]{LZ1} (also see \cite{S}).\\
(R1) The following limits exist finitely.
\begin{align*}
\alpha = \lim_{k \rightarrow \infty} \dfrac{b_k}{k}, \quad
\beta = \lim_{k \rightarrow \infty} \dfrac{a_k}{k}.	
\end{align*}
(R2) $\alpha = \beta^{-1}$ is a positive irrational number.	\\
(R3) Define a total order $>_{\alpha}$ on $\mathbb{Z}^2$ by setting
\begin{align*}
(i,j) >_{\alpha} (k,l) \iff i \alpha + j > k \alpha + l.	
\end{align*}
This order on $\mathbb{Z}^2$ is dense, which means that for every $(k,l) >_{\alpha} (0,0)$, there exist infinitely many $(i,j) \in \mathbb{Z}^2$ such that $(0,0) <_{\alpha} (i,j) <_{\alpha} (k,l)$.\\
Let us put $\mathbb{Z}^2(+) = \{(i,j) \in \mathbb{Z}^2 \ | \ (i,j) >_{\alpha} (0,0) \}$ and 
$\mathbb{Z}^2(-) = \{(i,j) \in \mathbb{Z}^2 \ | \ (i,j) <_{\alpha} (0,0) \}$.\\ 
(R4) $ \Lambda_{0} + (i,j) \in P(V) \implies \Lambda_{0} + (k,l) \in P(V) \ \forall \ (k,l) <_{\alpha} (i,j)$.\\
The ordering $>_{\alpha}$ naturally induces a triangular decomposition of $\tau$, say $\tau = \tau_{>_{\alpha}}^{-} \oplus \tau_{>_{\alpha}}^{0} \oplus \tau_{>_{\alpha}}^{+}$, where $\tau_{>_{\alpha}}^{0} = \mathfrak{g} \oplus \sum_{i=0}^{1} \mathbb{C}K_i \oplus \sum_{i=0}^{1} \mathbb{C}d_i$. \smallskip

\noindent \textbf{Assumption.} In our subsequent discussions (involving the case $n=1$) preceding Theorem \ref{T4.10}, we shall always assume that $V \ncong L(X)^A$  for any irreducible cuspidal $\tau_0$-module $X$ and $A \in GL(2, \mathbb{Z})$.

\blmma \label{Not equal}
$\tau_{-\underline{s}}v_{\mu} \neq (0)$ for any $\underline{s} = (a,b) \in \mathbb{Z}^2(+)$ and any non-zero $v_{\mu} \in V_{\mu}$.
\elmma
\begin{proof}
	From (R3), (R4) and (\ref{4.3}), it follows that for any $(i,j) \in \mathbb{Z}^2$, there exists $p \in \mathbb{Z}$ such that
	\begin{align}\label{(4.10)}
		\{m \in \mathbb{Z} \ | \ \Lambda_{0} + (i,j) + m(k,l) \in P(V) \} = (-\infty, p] \ \forall \ (k,l) \in \mathbb{Z}^2(+).
	\end{align}
	This implies that for any fixed $\underline{s} = (a,b) \in \mathbb{Z}^2(+)$, we have  $\tau_{r\underline{s}}v_{\mu} = (0)$ for $r \in
	\mathbb{N}$ large enough. \\
	If possible, let $\tau_{-\underline{s}}v_{\mu} = (0)$ for some non-zero vector $v_{\mu} \in V_{\mu}$. 
	Putting $c = \gcd(a,b) \in \mathbb{N}$, we have $\underline{s} =
	c(a^{\prime}, b^{\prime})$ where $\gcd(a^{\prime},b^{\prime}) = 1$. So
	there exist $p, q \in \mathbb{Z}$ such that $a^{\prime}q - b^{\prime}p =
	1$. Set $e_0^{\prime} = (a^{\prime},b^{\prime})$ and $e_1^{\prime} = (p,
	q)$, whence $\{e_0^{\prime}, \ e_1^{\prime}\}$ forms a $\mathbb{Z}$-basis of $\mathbb{Z}^2$. Now we can directly appeal to Lemma \ref{L4.6} to infer that, for any $0 \neq l \in \mathbb{Z}$, there exists $r_l \in \mathbb{Z}$ such that 
	\begin{align*}
	T_l := \{k \in \mathbb{Z} \ | \ \mu + le_1^{\prime} + ke_0^{\prime} \in P(V) \} = (-\infty, r_l] \ \text{or} \ [r_l, \infty).
	\end{align*}  
	It suffices to only consider $T_l = (-\infty, r_l]$ as a similar argument will also work for the other case. This implies that 
	$\tau_{le_1^{\prime} + cs_le_0^{\prime} \pm e_0^{\prime}}v_{\mu} =$ for large enough $s_l \in \mathbb{N}$. But as $\tau_{-ce_0^{\prime}}v_{\mu} = 0$, we can use the commutator relations on $\tau$ to deduce that  $\tau_{le_1^{\prime} \pm e_0^{\prime}}v_{\mu} = 0$ for any $0 \neq l \in \mathbb{Z}$. As a result, we have $\tau_{\pm (e_0^{\prime} + e_1^{\prime})}v_{\mu} = \tau_{\pm (e_0^{\prime} + 2e_1^{\prime})}v_{\mu} = 0$. It can be verified that $\{\tau_{\pm (e_0^{\prime} + e_1^{\prime})}, \tau_{\pm (e_0^{\prime} + 2e_1^{\prime})} \}$ generates $\tau$ as a Lie algebra, since $\{e_0^{\prime} + e_1^{\prime}, e_0^{\prime} + 2e_1^{\prime}\}$ is a $\mathbb{Z}$-basis of $\mathbb{Z}^2$. The irreducibility of $V$ then implies that $\tau$ acts trivially on $V$, which is a contradiction and hence the lemma is proved.   
\end{proof}

\blmma \label{L4.8}
$\big (\mu + \mathbb{Z}^{2}(+) \big) \bigcap P(V) \neq \emptyset \ \forall \ \mu \in P(V)$.
\elmma
\begin{proof}
	Let us assume the contrary. Then there exists $\mu \in P(V)$ with $\big (\mu +
	\mathbb{Z}^{2}(+) \big) \bigcap P(V) = \emptyset$ and thus $V$ is a highest weight module relative to the triangular decomposition in (R4) with its highest weight space $V_{+}^{> \alpha} = \{v \in V \ | \ \tau_{>_{\alpha}}^{+}v = 0 \} = (V_{+}^{> \alpha})_{\mu} \neq (0)$. Using the PBW theorem and the irreducibility of $V$, it can be shown that $V_{+}^{> \alpha}$ is an irreducible $\tau_{>_{\alpha}}^{0}$-module. 
	Now for any $k \in \mathbb{N}$, pick $\underline{m} = (a,b) \in \mathbb{Z}^2$ with $\dfrac{-1}{4k} <a \alpha + b < 0$. Due to (R3), there exist
	\begin{align}\label{(6.6)}
    \text{infinitely many} \ a \in \mathbb{Z} \ \text{and infinitely many} \ b \in \mathbb{Z} \ \text{with} \ (a,b) \ \text{satisfying the above relation}.
	\end{align}
	 Put $(a_i,b_i) = (2i-1)\underline{m} \ \forall \ 1 \leqslant i \leqslant k$. We can check that $(0,-1) = -e_1 <_{\alpha} (a_i,b_i) <_{\alpha} (0,0)$. Again set $\underline{p_i} = (1-2i)\underline{m} >_{\alpha} (0,0)$ and  $\underline{q_i} = e_1 - \underline{p_i} >_{\alpha} (0,0)$ for each $1 \leqslant i \leqslant k$. Observe that    
	\begin{align}\label{4.4}
		\underline{p_j} + \underline{q_j} = e_1 \ , \ \underline{q_1} >_{\alpha} \underline{q_i} \ , \ \underline{q_1} >_{\alpha}   \underline{p_j} \ \forall \ 2 \leqslant i \leqslant k, \ 1 \leqslant j \leqslant k, \\ 	
	\label{4.6}
		e_1 >_{\alpha} (0,0) , \	[\mathcal{Z},
		\mathcal{Z}] = (0) \ \text{and} \ \big (\mu +
		\mathbb{Z}^{2}(+) \big) \bigcap P(V) = \emptyset.
	\end{align}  
	\textbf{Claim 1.} $K_0$ and $K_1$ act trivially on $V_{+}^{> \alpha}$. \\
	If not, then let $K_1$ act by a non-zero scalar on $V_{+}^{> \alpha}$, say $\chi_1$. The proof of the other case involving $K_0$ is similar.
	For a fixed $0 \neq v \in V_{+}^{> \alpha}$ and any $k \in \mathbb{N}$,
	consider $\{(t^{-\underline{q_i}} K_1)(t^{-\underline{p_i}}
	K_1)v \}_{1 \leqslant i \leqslant k}$, all of which belong to
	$V_{\mu - e_1}$. Then these vectors cannot be linearly
	independent for each $k \in \mathbb{N}$, else we shall
	have $\dim V_{\mu -e_1} = \infty$. So there exists some $k \in
	\mathbb{N}$ satisfying
	$\sum_{i=1}^{k}\beta_i(t^{-\underline{q_i}}K_1)(t^{-\underline{p_i}}K_1)v
	= 0$ for some $\beta_1, \ldots, \beta_k \in \mathbb{C}$. We may assume that $\beta_1 \neq 0$. 
	Now set $\underline{r} = e_1 + 2\underline{m} >_{\alpha} \underline{0}$ and note that 
	\begin{align}\label{4.5}
	\underline{q_1} >_{\alpha} \underline{r} \ , \ \underline{q_i} <_{\alpha} \underline{r} \ , \ \underline{p_j} <_{\alpha} \underline{r} \ , \ e_1 >_{\alpha} \underline{r} \ \forall \ 2 \leqslant i \leqslant k, \ 1 \leqslant j \leqslant k.
	\end{align}
	An application of $t^{\underline{r}} d_0$ to the above equation
	together with (\ref{4.4}), (\ref{4.6}) and (\ref{4.5}) yields
	that 
	\begin{align*}
		a\beta_1(t^{\underline{r} - \underline{q_1}}K_1)(t^{-\underline{p_1}}K_1)v=0.
	\end{align*}       
	Again applying $t^{\underline{q_1} - \underline{r}}d_0$ to the above equation and using (\ref{4.4}), (\ref{4.6}) and (\ref{4.5}), we obtain
	\begin{align*}
		a^2\beta_1\chi_1(t^{-\underline{p_1}}K_1)v + a^2\beta_1(t^{\underline{r} - \underline{q_1}}K_1)(t^{\underline{q_1} -\underline{p_1} - \underline{r}}K_1)v=0.
	\end{align*}
	Applying $t^{\underline{p_1}}d_0$ to this equation, we obtain
	$\beta_1a^3\chi_{1}^2v = 0$. In view of (\ref{(6.6)}), we can now choose $a \neq 0$ and get $\beta_1 = 0$, which is a contradiction. Hence the claim. \\ 
	\noindent \textbf{Claim 2.} $d_0$ and $d_1$ act trivially on $V_{+}^{> \alpha}$. \\
	If not, then let $d_1$ act by the non-zero scalar $\mu(d_1)$ on $(V_{+}^{> \alpha})_{\mu}$. The other case involving $d_0$ can be dealt with analogously. Fix $0 \neq v \in V_{\mu}$ and consider $\{(t^{-\underline{q_i}}d_1)(t^{-\underline{p_i}}d_1)v \}_{1 \leqslant i \leqslant k} \subseteq V_{\mu - e_1}$ for any $k \in \mathbb{N}$. Then as in the last claim, there exist $k \in \mathbb{N}$ and $\beta_1, \ldots, \beta_k \in \mathbb{C}$ such that $\sum_{i=1}^{k} \beta_i (t^{-\underline{q_i}}d_1)(t^{-\underline{p_i}}d_1)v = 0$ for some non-zero $v \in V_{+}^{> \alpha}$, with $\beta_1 \neq 0$. Applying $t^{\underline{r}}d_1$ to the previous equation and using (\ref{4.4}), (\ref{4.6}), (\ref{4.5}), it can be deduced that there exists  $\beta_2^{\prime} \in \mathbb{C}$ such that 
	\begin{align}\label{(6.10)}
	\beta_1^{\prime}(t^{\underline{r} - \underline{q_1}}d_1)(t^{-\underline{p_1}}d_1)v + \beta_2^{\prime} (t^{\underline{r} - e_1}d_1)v + (*) = 0,  
	\end{align}
	where $\beta_1^{\prime} = (3b+2)\beta_1$ and $(*)$ consists of terms of the form $(t^{\underline{s}}K_i)(t^{\underline{k}}d_1)$ and $t^{\underline{s}}K_i$ ($ \underline{k}, \ \underline{s} \in \mathbb{Z}^{n+1}$ and $i = 0, 1$). Again applying $t^{e_1 - \underline{r}}d_1$ to (\ref{(6.10)}) and using (\ref{4.4}), (\ref{4.6}), (\ref{4.5}) along with Claim 1, we get 
	\begin{align}\label{(4.8)}
     3b^2\beta_1^{\prime} + 2b\beta_2^{\prime} = 0.
	\end{align}
Also applying $(t^{\underline{p_1}}d_1)(t^{\underline{q_1} - \underline{r}}d_1)$ to (\ref{(6.10)}) and using (\ref{4.4}), (\ref{4.6}), (\ref{4.5}) together with Claim 1 gives
\begin{align}\label{(6.12)}
2b^2(b + 2\mu(d_1))\beta_1^{\prime} + 3b^2\beta_2^{\prime} = 0.
\end{align}
Finally due to (\ref{(6.6)}), we can choose $b \in \mathbb{Z}$ such that the system of linear equations in (\ref{(4.8)}) and (\ref{(6.12)}) has a unique solution. This implies $\beta_1^{\prime} = \beta_2^{\prime} = 0$, which is a contradiction. Hence the claim. \smallskip

\noindent \textbf{Claim 3.} $\mathfrak{g}$ acts trivially on $V_{+}^{> \alpha}$. \\
If not, then there exist $v_0 \in V_{+}^{> \alpha}$ and $\lambda \in P_{\mathfrak{g}}^{+}$ such that $hv_0 = \lambda(h)v_0 \neq 0$ for some $h \in \mathfrak{h}$, where $\mathfrak{h}$ is a Cartan subalgebra of $\mathfrak{g}$. Then as in the previous claim, there exist some $k \in \mathbb{N}$ and $\beta_1, \ldots, \beta_k \in \mathbb{C}$ such that $\sum_{i=1}^{k} \beta_i (t^{-\underline{q_i}}d_1)(t^{-\underline{p_i}}d_1)v_0 = 0$, where $\beta_1 \neq 0$. Applying $h \otimes t^{\underline{r}}$ to the above equation and using (\ref{4.4}), (\ref{4.6}) and (\ref{4.5}), it can be deduced that there exists  $\beta_2^{\prime} \in \mathbb{C}$ such that
\begin{align} \label{(6.13)}
\beta_1^{\prime}(h \otimes t^{\underline{r} - \underline{q_1}})(t^{-\underline{p_1}}d_1)v_0 + \beta_2^{\prime}(h \otimes t^{\underline{r} - e_1})v_0 = 0, 
\end{align}   
where $\beta_1^{\prime} = (b+1)\beta_1$. Again applying $(h \otimes t^{\underline{p_1}})(t^{\underline{q_1} - \underline{r}}d_1)$ to (\ref{(6.13)}) and using (\ref{4.4}), (\ref{4.6}) and (\ref{4.5}) along with Claim 1 and Claim 2, we have $\beta_1^{\prime}b^2\lambda(h)^2v_0 = 0$. In view of (\ref{(6.6)}), we can now choose $b \neq 0$ and $b \neq -1$ to obtain $\beta_1 = \beta_1^{\prime} = 0$. This contradiction establishes our claim. 
Claims 1, 2 and 3 reveal that $\tau_{> \alpha}^{0}$ acts trivially on $V_{+}^{> \alpha}$. But since $V$ is an irreducible $\tau$-module,  $\tau$ must act trivially on $V$, which contradicts our initial assumption and so we are done.
\end{proof}

\bthm \label{T4.10}
$V \cong L(X)^A$ for some irreducible cuspidal module $X$ over $\tau_0$ and $A \in GL(2, \mathbb{Z})$.
\ethm	
\begin{proof}
	Let us assume the contrary. 
	Using Lemma \ref{L4.8}, pick $\underline{c} = (c_1,c_2) \in
	\mathbb{Z}^2(+)$ with $\Lambda_0 + \underline{c} \in P(V)$. By Lemma \ref{L3.8}, we must have $c_1c_2 \neq 0$ and thus there exists $s \in \mathbb{N}$ such that $(c_1,c_2) = s(c_1^{\prime},c_2^{\prime})$, where $c_1^{\prime}$ and $c_2^{\prime}$ are co-prime. So we can take $c_1$ and $c_2$ to be co-prime. 
	Now for any $k \in \mathbb{N}$, there
	exists $(a,b) \in \mathbb{Z}^2$ satisfying $0 < a\alpha + b <
	\dfrac{c_1\alpha + c_2}{4k}$, since $>_{\alpha}$ is a dense order
	on $\mathbb{Z}^2$ and $\underline{c} \in \mathbb{Z}^2(+)$. Due to (R3), it is clear that
	\begin{align}\label{(7.14)}
        \text{infinitely many} \ a \in \mathbb{Z} \ \text{and infinitely many} \ b \in \mathbb{Z} \ \text{with} \ (a,b) \ \text{satisfying the above relation}.
	\end{align} 
	Putting $\underline{m} = (a,b) \in \mathbb{Z}^2(+)$, set   
\begin{align*}
l = \max \{t \in \mathbb{Z} \ | \ \Lambda_0 + t \underline{c} \in P(V)
\}, \ l^{\prime} = \max \{t \in \mathbb{Z} \ | \ \Lambda_0 +
l\underline{c} + t \underline{m} \in P(V) \}. 
\end{align*}
	As $\Lambda_0 + \underline{c} \in P(V)$, it follows from (\ref{(4.10)}) that $l \in \mathbb{N}$ and $l^{\prime} \in \mathbb{Z}_{+}$.
	Let $\mu^{\prime} = \Lambda_0 + l\underline{c} + l^{\prime}\underline{m}$ and $(a_i^{\prime},b_i^{\prime}) = (1- 2i)\underline{m} \ \forall \ 1 \leqslant i \leqslant k$. We can now check that $-\underline{c} <_{\alpha} (a_i^{\prime},b_i^{\prime}) <_{\alpha} (0,0)$. Furthermore set $\underline{p_i}^{\prime} = (2i-1)\underline{m} >_{\alpha} (0,0)$ and  $\underline{q_i}^{\prime} = \underline{c} - \underline{p_i}^{\prime} >_{\alpha} (0,0)$ for each $1 \leqslant i \leqslant k$. 
	
	\noindent \textbf{Claim 1.} $\mu^{\prime} \neq 0$. \\
	Indeed, if $0 \in P(V)$, then we can take $\Lambda_{0} = 0$ and so the claim follows trivially as $(l\underline{c} + l^{\prime}\underline{m}) \in \mathbb{Z}^2(+)$ (see the paragraph before Theorem \ref{T4.2}). If $0 \notin P(V)$, then the claim is obvious as $\mu^{\prime} \in P(V)$. 
	
	\noindent \textbf{Claim 2.} $K_0$ and $K_1$ act trivially on $V_{\mu^{\prime}}$.\\
	If not, then let $K_1$ act non-trivially on $V_{\mu^{\prime}}$, say by $\chi_{1}$. A similar argument takes care of the case involving $K_0$. Now for every $k \in \mathbb{N}$ and any $0 \neq v \in V_{\mu^{\prime}}$, 
	$\{(t^{-\underline{q_i}^{\prime}}K_1)(t^{-\underline{p_i}^{\prime}}K_1)v \}_{1 \leqslant i \leqslant k} \subseteq V_{\mu^{\prime} - \underline{c}}$.
	Next we show that the collection $\{(t^{-\underline{q_i}^{\prime}}K_1)(t^{-\underline{p_i}^{\prime}}K_1)v \}_{1 \leqslant i \leqslant k}$ is not linearly independent for each $k \in \mathbb{N}$ and  $0 \neq v \in V_{\mu^{\prime}}$.
	If not, then $\dim V_{\mu^{\prime} - \underline{c}} \geqslant k$ for all $k \in \mathbb{N}$. Clearly
	$(l-1)\underline{c} + l^{\prime}\underline{m} \in \mathbb{Z}^2(+) \cup \{(0,0) \}$. If $(l-1)\underline{c} + l^{\prime}\underline{m} = 0$, then the assertion is obvious. On the other hand, if $(l-1)\underline{c} + l^{\prime}\underline{m} \in \mathbb{Z}^2(+)$, then we can apply Lemma \ref{Not equal} to infer that, for each $k \in \mathbb{N}$, there exists an injective linear map
	\begin{align*}
	\phi_k:  V_{\mu^{\prime} - \underline{c}} \longrightarrow (V_{\Lambda_0})^{\oplus N}, \ \text{where} \ N= \dim \tau_{(l-1)\underline{c} + l^{\prime}\underline{m}}.
	\end{align*}
	This implies that $N(\dim V_{\Lambda_0}) \geqslant k$ for every $k \in \mathbb{N}$, whence it follows that $\dim V_{\Lambda_0} = \infty$, which is a contradiction. This suggests that we can find some $k \in \mathbb{N}$ and $ \underline{0} \neq \underline{\beta} = (\beta_1, \ldots, \beta_k) \in \mathbb{C}^k$ satisfying
	$\sum_{i=1}^{k}\beta_i(t^{-\underline{q_i}^{\prime}}K_1)(t^{-\underline{p_i}^{\prime}}K_1)v  = 0$ for some $0 \neq v \in V_{\mu^{\prime}}$. Without loss of generality,  take $\beta_1 \neq 0$. Now setting $\underline{r}^{\prime} = \underline{c} - 2\underline{m}  >_{\alpha} \underline{0}$, note that      
	\begin{align}\label{(4.11)}
	[\mathcal{Z}, \mathcal{Z}] = (0), \ (\mu^{\prime} + \underline{m}) \notin P(V), \ (\mu^{\prime} + \underline{r}^{\prime} -\underline{p_j}^{\prime}) - (\mu^{\prime} + \underline{m}) \in \mathbb{Z}_2(+), \\ \label{(4.12)} (\mu^{\prime} + \underline{r}^{\prime} -\underline{q_i}^{\prime}) - (\mu^{\prime} + \underline{m}) \in \mathbb{Z}_2(+) \cup \{(0,0)\}  \ \forall \ 2 \leqslant i \leqslant k, \ 1 \leqslant j \leqslant k, \\ \label{(4.13)}
		(\mu^{\prime} + \underline{q_1}^{\prime} -\underline{r}^{\prime}) \notin P(V), \ (\mu^{\prime} + \underline{p_1}^{\prime}) \notin P(V), \ (\mu^{\prime} + \underline{c} -\underline{r}^{\prime}) \notin P(V).
	\end{align}
An application of $t^{\underline{r}^{\prime}} d_0$ to the above equation
together with (\ref{(4.11)}), (\ref{(4.12)}), (\ref{(4.13)}) and (R4) yields
\begin{align*}
(a-c_1)\beta_1(t^{\underline{r}^{\prime} - \underline{q_1}^{\prime}}K_1)(t^{-\underline{p_1}^{\prime}}K_1)v=0.
\end{align*}       
Again applying $t^{\underline{q_1}^{\prime} - \underline{r}^{\prime}}d_0$ to this equation and using (\ref{(4.11)}), (\ref{(4.12)}) and (\ref{(4.13)}), we obtain
\begin{align*}
a(a - c_1)\beta_1\chi_1(t^{-\underline{p_1}^{\prime}}K_1)v + a(a - c_1)\beta_1(t^{\underline{r}^{\prime} - \underline{q_1}^{\prime}}K_1)(t^{\underline{q_1}^{\prime} -\underline{p_1}^{\prime} - \underline{r}^{\prime}}K_1)v=0.
\end{align*}
Applying $t^{\underline{p_1}^{\prime}}d_0$ to this equation, we get
$\beta_1a^2(a - c_1)\chi_{1}^{2}v = 0$. By virtue of (\ref{(7.14)}), we can now choose $a \neq 0$ and $a \neq c_1$ to obtain $\beta_1 = 0$, which is a contradiction. Hence the claim. 

\noindent \textbf{Claim 3.} $d_0$ and $d_1$ act trivially on $V_{\mu^{\prime}}$. \\
For each $\underline{j} = (j_0,j_1) \in \mathbb{Z}^2 \setminus \{(0,0) \}$, set $d(\underline{j}) = j_0 t^{\underline{j}}d_1 - j_1 t^{\underline{j}}d_0.$ Then we have the following relations: 
\begin{align*}
[d(\underline{j}), d(\underline{k})] = (j_0k_1 - k_0j_1)d(\underline{j} + \underline{k}) + \eta(j_0k_1 - k_0j_1)^2 \sum_{s=0}^{1}j_s t^{\underline{j} + \underline{k}}K_s \ \text{for some} \ \eta 
\in \mathbb{C}; \\
[d(\underline{j}), t^{\underline{k}}K_i] = (j_0k_1 - k_0j_1)t^{\underline{j} + \underline{k}}K_i + (j_0 \delta_{i,1} - j_1 \delta_{i,0}) \sum_{s=0}^{1}j_s t^{\underline{j} + \underline{k}}K_s.
\end{align*} 
Now for any $k \in \mathbb{N}$ and $0 \neq v \in V_{\mu^{\prime}}$, consider $\{d(-{\underline{q_i}^{\prime}}) d({-\underline{p_i}^{\prime}})v \}_{1 \leqslant i \leqslant k} \subseteq V_{\mu^{\prime} - \underline{c}}$. Then we can again invoke the aforesaid argument to prove that these vectors cannot be linearly independent for each $k \in \mathbb{N}$ and non-zero $v \in V_{\mu^{\prime}}$. So there exists $ \underline{0} \neq \underline{\beta} = (\beta_1, \ldots, \beta_k) \in \mathbb{C}^k$ with $k \in \mathbb{N}$ such that
	$\sum_{i=1}^{k}\beta_i \big(d(-{\underline{q_i}^{\prime}})  d({-\underline{p_i}^{\prime}}) \big)v = 0$ for some $v \in V_{\mu^{\prime}}$. Without loss of generality, let $\beta_1 \neq 0$. Applying $d({\underline{r}^{\prime}})$ to the previous equation and using (\ref{(4.11)}), (\ref{(4.12)}), (\ref{(4.13)}) and (R4), it can be deduced that there exist  $\eta_0, \eta_1, \eta_0^{\prime}, \eta_1^{\prime} \in \mathbb{C}$ such that 
	\begin{align}\label{(7.10)}
    \beta_1 \big ((ac_2-bc_1)d(-\underline{m}) + \sum_{i=0}^{1}\eta_i (t^{-\underline{m}}K_i) \big) d(-\underline{m})v + \beta_1^{\prime} \big (\sum_{i=0}^{1}\eta_i^{\prime} (t^{-2 \underline{m}}K_i) \big)v = 0,  
	\end{align}
	where we take $\beta_1^{\prime} = \sum_{i=2}^{k} \beta_i$. Again apply $t^{2\underline{m}}d_1$ to (\ref{(7.10)}) and use (\ref{(4.11)}), (\ref{(4.12)}), (\ref{(4.13)}) along with  Claim 2 to deduce that 
	\begin{align}\label{(7.8)}
	\beta_1ab(ac_2-bc_1) \big ((2b^2 + b + 2)\mu^{\prime}(d_0) + a(2b + 5)\mu^{\prime}(d_1) \big )v = 0.
	\end{align}
	Next we show that $ac_2 \neq bc_1$. If not, then $ac_2=bc_1$. But as $c_1$ and $c_2$ are co-prime, we get $(a,b) = s^{\prime}(c_1,c_2)$ for some $s^{\prime} \in \mathbb{Z}$, which clearly contradicts our choice of $\underline{m}$. So we must have $ac_2 \neq bc_1$, which thereby reduces (\ref{(7.8)}) to
	\begin{align}\label{(7.9)}
	ab \big((2b^2 + b + 2)\mu^{\prime}(d_0) + a(2b + 5)\mu^{\prime}(d_1) \big) = 0
	\end{align}
	An application of $(t^{\underline{m}}d_1)(t^{\underline{m}}d_1)$ to (\ref{(7.10)}) together with (\ref{(4.11)}), (\ref{(4.12)}), (\ref{(4.13)}) and Claim 2 yields
	\begin{align}\label{(7.12)}
	2\beta_1b^2(ac_2-bc_1)(a \mu^{\prime}(d_1) - b\mu^{\prime}(d_0))^2 = 0.
	\end{align}
	Now (\ref{(7.14)}) allows us to pick $a, b \in \mathbb{Z}$ such that (\ref{(7.9)}) and (\ref{(7.12)}) are finally reduced to 
	\begin{align}\label{(7.22)}
	(2b^2 + b + 2)\mu^{\prime}(d_0) + a(2b + 5)\mu^{\prime}(d_1) = 0 = b \mu^{\prime}(d_0) - a\mu^{\prime}(d_1). 
	\end{align}
	Then we can again apply (\ref{(7.14)}) to appropriately choose $a,b \in \mathbb{Z}$ such that the system of linear equations in (\ref{(7.22)}) has a unique solution $\mu^{\prime}(d_0) = 0 = \mu^{\prime}(d_1)$. Hence the claim. 
	
	But this contradicts Claim 1, which henceforth proves the theorem. 

\end{proof} 

Our Theorem \ref{T4.2} can be now deduced in a more or less similar manner as in \cite[Lemma 3.8]{LZ1} and \cite[Theorem 3.9]{LZ1} by suitably applying Lemma \ref{L2.2}, Lemma \ref{L3.8}, Proposition \ref{P4.3}, Theorem \ref{T4.10} at various stages, together with Remark \ref{Dense}(1),  Remark \ref{R5.3} and then performing appropriate change of co-ordinates. But for the sake of completeness, we just outline the main steps of the proof while referring to \cite[Lemma 3.8]{LZ1} and \cite[Theorem 3.9]{LZ1} for the complete details. \\ 

\noindent \textbf{Proof of Theorem \ref{T4.2}.} We shall prove this theorem by induction on $n$. Note that for $n=1$, this is precisely Theorem \ref{T4.10}. Now suppose that the theorem holds good for any $n \leqslant N-1$, where $N \geqslant 2$. We shall prove this theorem for $n=N$.\\
\textbf{Step 1.} If $L(X)$ is a non-trivial Harish-Chandra module of the highest weight type (see Section \ref{Highest Weight}), then there exists $\lambda \in D^*$ such that $P(L(X)) = P(X) \bigcup (\lambda - \mathbb{N}e_0 + \mathbb{Z}^n)$, where $\mathbb{Z}^n \cong \text{span} \{e_1, \ldots, e_n \}$. \\
$Proof \ of \ Step \ 1.$ Assume the contrary and then proceed as in 
\cite[Lemma 3.8]{LZ1}. Replace the rank $2$ Virasoro algebra
by the corresponding full toroidal Lie algebra of rank $2$ and then invoke Lemma \ref{L2.2} to obtain a non-trivial irreducible subquotient, which is a non-trivial GHW module, thanks to Remark \ref{Dense}. Finally using our Lemma \ref{L3.8}, we obtain a contradiction, similar to \cite[Lemma 3.8]{LZ1}. \\
\textbf{Step 2.} If $(\Lambda_0 + \underline{m} + G_0) \cap P(V) \subseteq \{0 \}$ for some $\underline{m} \in \mathbb{Z}^{n+1}$ and a co-rank $1$ subgroup $G_0$ of $\mathbb{Z}^{n+1}$, then the theorem holds good.\\         
$Proof \ of \ Step \ 2.$ This step follows using our Proposition \ref{P4.3} (refer to the paragraph above (3.20) of \cite[Theorem 3.9]{LZ1} for more details).

By Step 2, we can assume that $(\Lambda_0 + \underline{m} + G_0) \cap P(V) \nsubseteq \{0 \}$ for any $\underline{m} \in \mathbb{Z}^{n+1}$ and any co-rank $1$ subgroup $G_0$ of $\mathbb{Z}^{n+1}$.
Consider $V_{\Lambda_0 + \underline{m} + G_0} = \bigoplus_{g \in G_0} V_{\Lambda_0 + \underline{m} + g}$, which is a Harish-Chandra module over a full toroidal Lie algebra of rank $(N-1)$ (instead of Vir$[G_0$] in \cite[Theorem 3.9]{LZ1}) and then apply Lemma \ref{L2.2} and our inductive hypothesis along with Step 1 to conclude that for any $\underline{m} \in \mathbb{Z}^{n+1}$ and any co-rank $1$ subgroup $G_0$ of $\mathbb{Z}^{n+1}$, there exist a subgroup $G_{0,1}$ of $G_0$, $\lambda_0^{\prime} \in \Lambda_0 + \underline{m} + G_0$ and $g_{0,1} \in G_0 \setminus \{0 \}$ with $G_0 = \mathbb{Z}g_{0,1} \oplus G_{0,1}$ such that   $\lambda_0^{\prime} + G_{0,1} - \mathbb{N}g_{0,1} \subseteq P(V)$.

\noindent \textbf{Step 3.} There exist no $\lambda_0 \in P(V), \ t_0 \in \mathbb{Z}, \ \underline{m_0}, \underline{m_1} \in \mathbb{Z}^{n+1} \setminus \{\underline{0} \}$ or subgroups $G_1^{\prime} \subseteq G_0^{\prime} \subseteq \mathbb{Z}^{n+1}$ with $\mathbb{Z}^{n+1} = \mathbb{Z} \underline{m_0} \oplus G_0^{\prime}$ and $G_0^{\prime} = \mathbb{Z} \underline{m_1} \oplus G_1^{\prime}$  satisfying
\begin{align*}
\lambda_0 - \mathbb{Z}_{+} \underline{m_1} + G_1^{\prime}, \ \lambda_0 + t_0 \underline{m_1} + \mathbb{Z}_{+}\underline{m_1} +  G_1^{\prime} \subseteq P(V) \ (\text{if} \ t_0 \leqslant0, \ \text{then} \ \lambda_0 +  G_0^{\prime} \subseteq P(V)).
\end{align*}
$Proof \ of \ Step \ 3.$ This assertion can be proved by essentially appealing to the same arguments presented in Claim 1 of \cite[Theorem 3.9]{LZ1}.

For each $t \in \mathbb{Z}$, define $\overline{G}_t := te_0 + \mathbb{Z}e_1 + \mathbb{Z}e_2 + \ldots + \mathbb{Z}e_N$.\\  
\textbf{Step 4.} If for $\lambda_0 \in \Lambda_0 + \mathbb{Z}^{n+1}, \underline{m_1}, \underline{m_1}^{\prime} \in \overline{G}_0 \setminus \{0\}$ and subgroups $G_1$ and $G_1^{\prime}$ of $\overline{G}_0$,  with $\overline{G}_0 = \mathbb{Z} \underline{m_1} \oplus G_1, \ \overline{G}_0 = \mathbb{Z}\underline{m_1}^{\prime} \oplus G_1^{\prime}$, we have 
\begin{align*}
\lambda_0 - \mathbb{N} \underline{m_1} + G_1, \ \lambda_0 - \mathbb{N} \underline{m_1}^{\prime} + G_1^{\prime} \subseteq P(V), \ \text{then} \ G_1 = G_1^{\prime}.
\end{align*}   
$Proof \ of \ Step \ 4.$ The proof is similar to Claim 2 of \cite[Theorem 3.9]{LZ1}.

Set $V_{\Lambda_0 + \overline{G}_t} = \bigoplus _{\underline{k} \in \overline{G}_t} V_{\Lambda_0 + \underline{k}}$ for each $t \in \mathbb{Z}$, which is a Harish-Chandra module over a full toroidal Lie algebra of rank $(N-1)$.\\
\textbf{Step 5.} There exist a co-rank 1 subgroup $G_0$ in $\overline{G}_0, \ \alpha_t \in \Lambda_0 + \overline{G}_t$ and $\underline{k_0} \in \overline{G}_0$ with $\overline{G}_0 = \mathbb{Z}_{+}\underline{k_0} \oplus G_0$ such that either
\begin{align*}
P(V_{\Lambda_0 + \overline{G}_t}) \setminus \{0 \} = &\ (\alpha_t +
\mathbb{Z}_{+}\underline{k_0} + G_0) \setminus \{0 \} \ \text{or} \
(\alpha_t -
\mathbb{Z}_{+}\underline{k_0} + G_0) \setminus \{0 \}.
\end{align*}  
$Proof \ of \ Step \ 5.$ The proof is precisely as in (3.42) of \cite[Theorem 3.9]{LZ1}.

Now proceed similarly as in \cite[Theorem 3.9]{LZ1} to finally get a contradiction to Step 3, which henceforth proves the theorem. 


\section{The Main Classification Theorems} \label{Theorems}
In this section, we use results from Section \ref{Central}, Section \ref{GHW modules}, Section \ref{Highest Weight} and Section \ref{S6} to finally classify the irreducible Harish-Chandra modules over the full toroidal Lie algebra. As one of the direct applications, we also obtain the classification of all possible irreducible Harish-Chandra modules over the higher-dimensional Virasoro algebra, thereby proving Eswara Rao's conjecture \cite{E}.

\bthm \label{C1}
Let $V$ be a non-trivial irreducible Harish-Chandra module over
$\tau$. Then:
\begin{enumerate}
\item $V$ is either a cuspidal module or a highest weight type module.
\item If $V$ is cuspidal, then $V$ is isomorphic to one of the following.\\
(a) $V(c, \lambda_1, \lambda_2, \underline{\alpha})$ for some quadruplet $(c, \lambda_1, \lambda_2, \underline{\alpha}) \in \mathbb{C} \times (P_{\mathfrak{g}}^{+})^{\times} \times P_{\mathfrak{sl}_{n+1}}^{+} \times \mathbb{C}^{n+1}$, where $(P_{\mathfrak{g}}^{+})^{\times} = (P_{\mathfrak{g}}^{+}) \setminus \{0\}$.\\
(b) A $\mathcal{W}_{n+1}$-module $V(c, 0, \lambda_2, \underline{\alpha})$, with $(c, \lambda_2, \underline{\alpha}) \in \mathbb{C} \times P_{\mathfrak{sl}_{n+1}}^{+} \times \mathbb{C}^{n+1}$, where $(\lambda_2, c) \neq (\omega_k,k)$ for any fundamental weight $\omega_k$ of $\mathfrak{sl}_{n+1}$ and $1 \leqslant k \leqslant n$.\\ 
(c) A $\mathcal{W}_{n+1}$-module $V(n+1, 0, 0, \underline{\alpha})$, where $\underline{\alpha} \notin \mathbb{Z}^{n+1}$.\\ 
(d) A $\mathcal{W}_{n+1}$-submodule $\pi_k \big( V(k , 0, \omega_k, \underline{\alpha}) \big)$ of $V(k+1 , 0, \omega_{k+1}, \underline{\alpha})$ for some $\underline{\alpha} \in \mathbb{C}^{n+1}$, where $0 \leqslant k \leqslant n$ and $\omega_0 = \omega_{n+1} = 0$. 
\item A highest weight type module is isomorphic to $L(V^{\prime}(a, b, c, \lambda_1, \lambda_2, \underline{\alpha}))^A$ for $(a,b,c, \lambda_1, \lambda_2, \underline{\alpha}) \in \mathbb{C} \times \mathbb{C} \times \mathbb{C} \times P_{\mathfrak{g}}^{+} \times P_{\mathfrak{sl}_n}^{+} \times \mathbb{C}^{n}$ and $A \in GL(n+1, \mathbb{Z})$.   
\end{enumerate}

\begin{proof}
(1) follows from Theorem \ref{T2.9} and Theorem \ref{T4.2}, whereas (2) follows from Theorem \ref{T2.14}.\\
(3) follows immediately from Theorem \ref{L(X)} and Remark \ref{R5.3}. 
\end{proof}
\ethm
\brmk \label{Bounded}
 Theorem \ref{C1} unravels that the irreducible \textit{bounded} modules (see Remark \ref{Yuly}) are all isomorphic to $L(V^{\prime}(a, b, c, \lambda_1, \lambda_2, \underline{\alpha}))$ for some $(a,b,c, \lambda_1, \lambda_2, \underline{\alpha}) \in \mathbb{C} \times \mathbb{C} \times \mathbb{C} \times P_{\mathfrak{g}}^{+} \times P_{\mathfrak{sl}_n}^{+} \times \mathbb{C}^{n}$. This gives all the simple objects in the category of bounded modules, which was introduced in \cite{B}. 
\ermk

\bthm \label{C2}
Let $V$ be a non-trivial irreducible Harish-Chandra module over
$\mathcal{V}ir$. Then:
\begin{enumerate}
	\item $V$ is either a cuspidal module or a highest weight type module.
	\item If $V$ is cuspidal, then $V$ is isomorphic to one of the following.\\
	(a) A $\mathcal{W}_{n+1}$-module $V(c, \lambda_2, \underline{\alpha})$, with $(c, \lambda_2, \underline{\alpha}) \in \mathbb{C} \times P_{\mathfrak{sl}_{n+1}}^{+} \times \mathbb{C}^{n+1}$, where $(\lambda_2, c) \neq (\omega_k,k)$ for any fundamental weight $\omega_k$ of $\mathfrak{sl}_{n+1}$ and $1 \leqslant k \leqslant n$.\\ 
	(b) A $\mathcal{W}_{n+1}$-module $V(n+1, 0, \underline{\alpha})$, where $\underline{\alpha} \notin \mathbb{Z}^{n+1}$.\\ 
	(c) A $\mathcal{W}_{n+1}$-submodule $\pi_k \big(V(k, \omega_k, \underline{\alpha}) \big)$ of $V(k+1, \omega_{k+1}, \underline{\alpha})$, where $\underline{\alpha} \in \mathbb{C}^{n+1}$, $0 \leqslant k \leqslant n$ and $\omega_0 = \omega_{n+1} = 0$. 
	\item A highest weight type module is isomorphic to $L(V^{\prime}(a, b, c, \lambda, \underline{\alpha}))^A$ for some $(a,b,c, \lambda, \underline{\alpha}) \in \mathbb{C} \times \mathbb{C} \times \mathbb{C} \times P_{\mathfrak{sl}_n}^{+} \times \mathbb{C}^{n}$ and $A \in GL(n+1, \mathbb{Z})$.   
\end{enumerate}
\ethm
\begin{proof}
	(1) follows from Theorem \ref{T2.9} and Theorem \ref{T4.2} by taking $\mathfrak{g} = (0)$, whereas (2) follows from Theorem \ref{T2.13} and \cite[Theorem 5.4]{BF}. (3) follows similarly from Theorem \ref{L(X)} and Remark \ref{R5.3}.  
\end{proof}
\brmk \label{Analogue}
Our techniques can be directly utilized to also classify the irreducible
Harish-Chandra modules over $\mathcal{H} \mathcal{V}ir = \mathfrak{h}
\otimes A \oplus \mathcal{Z} \oplus \mathcal{W}_{n+1}$, where
$\mathfrak{h}$ is a Cartan subalgebra of $\mathfrak{g}$, by 
replacing the finite-dimensional irreducible $\mathfrak{g}$-module by a
$1$-dimensional $\mathfrak{h}$-module. This Lie algebra can be viewed as a higher-dimensional generalization of the twisted Heisenberg--Virasoro algebra admitting a common extension (obtained by taking $C_D = C_I$ and $C_{DI}=0$ in the definition of HVir in Subsection \ref{Twisted}). 
\ermk

\brmk
While writing this paper, we found out that the Harish-Chandra modules over $\mathcal{V}ir$ were also studied in \cite{JJ}. However, the techniques used in that paper are largely limited to the framework of \textit{non-zero} level modules, due to which the authors were able to classify only those irreducible Harish-Chandra modules, where at least one of the $K_i$'s act \textit{non-trivially}. In the current paper, we take a completely different approach to that of \cite{JJ}, which helps us to provide a uniform proof that works not only for the higher-dimensional Virasoro algebra but also for the full toroidal Lie algebra, in case of \textit{both level zero as well as non-zero level modules}.  
\ermk

\smallskip

\noindent \textbf{Acknowledgements.} The author would like to thank Prof. Eswara Rao for suggesting the problem and for some fruitful discussions. He also thanks Prof. Apoorva Khare for going through the whole paper and for his helpful suggestions. Most of this work was done when the author was a Research Associate at the Theoretical Statistics and Mathematics Unit, Indian Statistical Institute, Bangalore. The author gratefully acknowledges the financial support and the excellent working conditions provided by the institute.

\end{document}